\documentclass[a4paper,11pt,leqno]{amsart}
 
\usepackage{amssymb,hyperref}%,showkeys}
\theoremstyle{plain}

\newtheorem{theo}{Theorem}[section]

\newtheorem{proposition}[theo]{Proposition}
\newtheorem{corollaire}[theo]{Corollary}

\newtheorem{conjecture}[theo]{Conjecture}
\newtheorem{theorem}[theo]{Theorem}
\newtheorem{lemma}[theo]{Lemma}
\newtheorem{remark}[theo]{Remark}

\theoremstyle{definition}

\theoremstyle{remark}

 \def\RR{{\mathbb R}}
\def\CC{{\mathbb C}}
\def\ZZ{{\mathbb Z}}

\def\NN{{\mathbb N}}

\begin{document}

\title[Holomorphic metrics]{Global rigidity of holomorphic Riemannian metrics on compact complex 3-manifolds}

\author[S. Dumitrescu, A. Zeghib]{Sorin DUMITRESCU$^\star$ \& Abdelghani ZEGHIB$^\dagger$}

\address{${}^\star$ D\'epartement de Math\'ematiques d'Orsay, 
\'equipe de Topologie et Dynamique,
Bat. 425, U.M.R.   8628  C.N.R.S.,
Univ. Paris-Sud (11),
91405 Orsay Cedex, France}
\email{Sorin.Dumitrescu@math.u-psud.fr}
 
\address{${}^\dagger$ CNRS, UMPA, \'ecole
Normale Sup\'erieure de Lyon, France}
\email{zeghib@umpa.ens-lyon.fr}

\thanks{}
\date{\today}

\setcounter{tocdepth}{1}

\maketitle

\begin{abstract}    We study compact complex 3-manifolds  admitting holomorphic Riemannian metrics.
We prove a uniformization result:   up to  a finite   unramified cover, such a manifold   admits  a holomorphic Riemannian metric   of
  constant sectionnal curvature.  \\ 
\end{abstract}

\section{Introduction}

A {\it holomorphic Riemannian metric} $g$ on a complex manifold $M$ is a holomorphic field of non degenerate  complex quadratic forms on the holomorphic tangent bundle $TM$.
Formally,  $g$ is a holomorphic section of the bundle $S^2(T^{*}M)$ such that $g(m)$ is non degenerate for all $m \in M$. This has nothing to do with the more usual Hermitian metrics.
It is in fact nothing but the complex version of Riemannian metrics. 
Observe that  since complex quadratic forms  have no signature, there is here no distinction   between the Riemannian and pseudo-Riemannian cases. This  observation was the origin of the nice use by F. Gau{\ss}  of the complexification technic of (analytic)  Riemannian metrics on surfaces, in order to find  conformal coordinates for them. Actually, the complexification of analytic Riemannian metrics leading to holomorphic ones, is becoming a standard trick (see for instance~\cite{Fra}).

As in the real case,  a  holomorphic Riemannian metric on $M$ gives rise to  a covariant differential calculus, i.e.    a  Levi-Civita (holomorphic) linear
connection, and to geometric features:   curvature tensors, geodesic (complex) curves~\cite{Le, Leb}.

Locally,  a holomorphic Riemannian metric has the form 
$\Sigma g_{ij}(z)dz_i dz_j$, where $(g_{ij}(z))$ is a complex inversible   symmetric matrix depending holomorphically on $z$.
The standard   example is that of  the global flat  holomorphic Riemannian metric  $dz_{1}^2+dz_{2}^2+ \ldots +dz_{n}^2$ on $\CC^n$. This metric is translation-invariant and thus 
goes down to any quotient of $\CC^n$ by a lattice. Hence complex torii  possess  (flat) holomorphic Riemannian metrics. This is however a very  special situation since, contrary to 
 real case, only  few   {\it compact} complex manifolds admit holomorphic Riemannian metrics. 
Our  goal in this paper is to illustrate this rigidity by the following uniformization   theorem:

\begin{theorem}  \label{principal}
 If a compact connected complex  3-manifold  $M$ admits a holomorphic Riemannian metric, then, up to  a finite   unramified cover,  $M$   admits a holomorphic Riemannian
 metric of constant sectionnal curvature. 
\end{theorem}

The starting   point of this result is the main result  of~\cite{Dum1}: 

\begin{theorem}  \label{Dumitrescu}  \cite{Dum1}
Any holomorphic Riemannian metric on a compact connected complex $3$-manifold is locally homogeneous. More generally, 
 if  a compact connected complex $3$-manifold $M$ admits a   holomorphic Riemannian metric, then any  holomorphic geometric structure of affine type on $M$
is locally homogeneous. 
\end{theorem}

 The simplest complex compact manifolds endowed with  holomorphic Riemannian metrics are those obtained as 
  a (left) quotient of a complex Lie group $G$ by a co-compact lattice $\Gamma$. The holomorphic Riemannian metric on $G$ is left invariant and  can be constructed  by left translating 
  any complex non-degenerate quadratic form defined on the Lie algebra ${\mathcal G}$.  For such (special)  spaces, our result follows from the following 
  ``algebraic''  fact:

\begin{proposition}  \label{unimodulaire} A 3-dimensional unimodular complex Lie group 
admits a left invariant holomorphic Riemannian metric of constant sectionnal curvature. This metric is flat exactly when the group is solvable. 
\end{proposition}

\begin{remark}
This   is just  the complexified version of the fact that any real unimodular $3$-dimensional Lie group admits a left invariant  pseudo-Riemannian metric (which is thus either Riemannian or  Lorentzian)  of 
constant sectionnal curvature \cite{Mi}, \cite{Ra} (see also \S \ref{Examples}).
\end{remark}

The main  result of this paper can be seen as a generalization of the previous proposition. 
More precisely, we prove:

\begin{theorem} \label{result}
Let $M$ be a compact connected complex $3$-manifold which admits a (locally homogeneous) holomorphic Riemannian metric $g$. Then: 

(i) If the Killing Lie algebra of $g$ has a non trivial  semi-simple part, then it  preserves some  holomorphic Riemannian metric on $M$ with  constant sectionnal curvature.

(ii) If the Killing Lie algebra of $g$ is solvable, then, 
up to  a finite  unramified cover,
 $M$   is a quotient $\Gamma \backslash G$, where 
$\Gamma$ is a lattice of $G$ and $G$ is either  the complex Heisenberg group,  or the complex $SOL$ group.
Furthermore, the pull-back of $g$ on the universal cover  of $M$ is a left invariant holomorphic Riemannian metric on $G$.
\end{theorem}

Note that the group $SOL$ is the complexification of  the  affine isometry group of the Minkowski plane $\RR^{1, 1}$ or equivalently the isometry group of $\CC^2$ endowed with its flat holomorphic Riemannian  metric (see \S \ref{Examples}).

\subsection{Completeness} \label{question1}
Our present result  does not end the story, essentially because of remaining {\it completeness}  questions, and those on the algebraic structure of the fundamental group. 

{\it  It remains to classify the compact complex $3$-manifolds endowed with a holomorphic Riemannian metric of  constant sectionnal curvature}.  

Let us give details in the flat case. So, let $M$ be a compact manifold locally modelled on the flat model ${\mathbb C}^3$. With Thurston's  terminology~\cite{Th}, $M$ admits 
a $(O(3, {\mathbb C}) \ltimes {\mathbb C}^3, {\mathbb C}^3)$-structure.   The challenge  remains: \\

1) {\it Markus conjecture:}  Is $M$ complete, i.e. is there  $\Gamma \subset O(3, {\mathbb C}) \ltimes {\mathbb C}^3$ acting properly discontinuously on ${\mathbb C}^3$
such that $M = {\mathbb C}^3 / \Gamma$ ? (see \cite{Mar}).\\

2) {\it Auslander conjecture: }  Assuming $M$ as above, is $\Gamma$ solvable? 

Note that these questions are settled in the setting of (real) flat Lorentz manifolds~\cite{Car, Fried-Goldman}, but unsolved  for general (real) pseudo-Riemannian metrics. The real part
of the holomorphic Riemannian metric is a (real) pseudo-Riemannian metric of signature $(3,3)$ for which both previous conjectures are still open.

More details  about completeness in the case of a non-zero constant sectionnal curvature are  in  \S  \ref{non-zero curvature}.

\subsubsection*{Comparison with \cite{DZ}} The present article is naturally linked to   
  our recent work on  the classification of essential lorentz geometries in dimension $3$. 
   There are  similarities in the
algebraic classification of all  possible local Killing algebras.
 However, we had to modify significantly    our methods because in~\cite{DZ}  we used global results about the  classification of (real)  Riemannian
Killing fields \cite{Ca} and about the classification of non-equicontinuous Lorentz Killing fields \cite{Zeg} which do not  exist in the holomorphic setting.

\subsubsection*{Related works}  There are various  works dealing with different holomorphic geometric structures, and sharing the same philosophy as ours here, that is, a ``strong global  rigidity'' of such objects on compact complex manifolds. As an example, we can quote 
\cite{IKO, KO, HM, D1}, and especially~\cite{Radloff}, about  holomorphic conformal structures on projective 3-manifolds.  As an extension of both their results and ours, we believe   a global rigidity result is true  for  holomorphic conformal metrics in the  framework   of complex (not necessarily projective)  3-manifolds.

 \subsection{Plan of the proof}  We   briefly indicate the important steps in the proof of Theorem \ref{result}. Thanks to theorem~\ref{Dumitrescu}, we are in a locally homogeneous situation: our manifold $M$ is locally modelled on a $(G, G/I)$-geometry  in Thurston's sense~\cite{Th}, where $I$ is a closed subgroup of the Lie group $G$. We  have two objects to understand:
 \begin{enumerate}
 
 \item $G$ and $I$ inside it;

 \item the holomony morphism $\rho: \pi_1(M) \to G$.
  \end{enumerate}
  
 The first step consists on finding all 3-dimensional complex homogeneous spaces $G/I$
 such that  the $G$-action on  $G/I$ preserves some holomorphic Riemannian metric (i.e. the adjoint representation of $I$ preserves some non-degenerate
 complex quadratic form on the quotient $\mathcal G /  \mathcal I$ of the corresponding Lie algebras). Despite a ``quick'' reduction to the case where 
 $G$ has dimension 4 and is solvable, our  solution needs a geometric tool which is the existence of a codimension one   geodesic foliation $\mathcal F$.
 
 The second part is a standard problem: classify compact manifolds locally modelled  on  a given $(G, G/I)$-geometry. It has two sides. The first is  completeness, that is the holonomy group $\Gamma= \rho (\pi_1(M))$ acts properly on $G/I$  and $M$ is a compact quotient  $\Gamma  \backslash G /I$. The second side  classify the discrete groups  $\Gamma$. If $G$ is solvable, we prove that $M$ is complete and, up to a finite cover, it is a quotient of $Heis$ or $SOL$ by a lattice.

\section{Examples}    \label{Examples}

A first obstruction to the existence of a holomorphic Riemannian metric on a compact complex manifold is  its  first Chern class. Indeed, a holomorphic Riemannian
metric on $M$ provides an isomorphism between $TM$ and $T^{*}M$. In particular, the canonical bundle $K$ is isomorphic to the anti-canonical bundle $K^{-1}$
and $K^2$ is trivial. This means that the first Chern class of $M$ vanishes and, up to a double unramified cover, $M$ possesses a holomorphic volume form.

\subsection*{Quadratic differentials}  

The previous obstruction implies that the only Riemann surfaces (complex curves)  which admit  (1-dimensional) holomorphic Riemannian metrics are elliptic curves.

A  (holomorphic) quadratic differential on a Riemann surface has locally the form $\phi(z)dz^2$, where $\phi$ is a holomorphic function.  It can be seen as a ``singular'' holomorphic Riemannian 
metric. Outside its null set, it determines  a 
 holomorphic Riemannian metric  which 
 is flat, i.e. locally isomorphic to $dz^2$ (similarly to the situation of real 1-dimensional Riemannian metrics). 
 This also endows the surface with    a translation-structure (i.e. a $(\CC, \CC)$-structure in Thurston's sense), which are nowadays a central subject of study from various points of view (see for instance \cite{Strebel, KZ, McMullen}...).
 
 In  higher dimension, a quadratic differential can be defined as a holomorphic  section of $S^2(T^*M)$. However  no systematic study of them seems to exist, even in the case of surfaces or 3-manifolds. One motivation of our interest to  holomorphic Riemannian metrics, is that they  correspond exactly to the 
 case where this quadratic differential is  non-degenerate. This is surely 
 a strong hypothesis, but our rigidity results give evidence that  other more flexible cases can also be handled.

 \subsection*{Kaehler case}

We have seen above  that complex torii  admit flat holomorphic Riemannian metrics. In fact, up to an unramified finite cover, they are the only compact Kaehler manifolds
admitting holomorphic Riemannian metrics~\cite{IKO}. 

\subsection*{Surface case} 

 In the surface case (Kaehler or not), the sectionnal curvature is a holomorphic {\it function}, and thus constant by compactness. It was proved in~\cite{Dum} that this curvature must in fact vanish and, up to an unramified finite cover, only   complex torii admit  holomorphic Riemannian metric. In particular, there is no compact surface having a holomorphic Riemannian
 metric of non zero constant sectionnal curvature.

 \subsection*{Universal holomorphic  Riemannian  spaces of constant curvature} One can multiply a holomorphic Riemannian metric by a complex constant $\lambda$ which induces a multiplication 
 by $\lambda^{-2}$ of its sectionnal curvature. Therefore, only 
 the vanishing or not (but not the sign)  of the curvature
 is relevant.

 \subsubsection*{The flat case}  The model   $(\CC^n, dz_{1}^2+dz_{2}^2+ \ldots +dz_{n}^2)$ is (up to isometry) the unique    $n$-dimensional complex simply-connected manifold endowed with a flat and geodesically complete holomorphic Riemannian metric. Its isometry group   is $O(n, \CC) \ltimes \CC^n$.  Any flat
 holomorphic Riemannian metric on a complex manifold of dimension $n$ is locally  isometric to this model, equivalently, it has a 
  $(O(n, \CC) \ltimes \CC^n,   \CC^n)$-structure~\cite{Th,Wo}. 
 This geometry  can be seen as a complexification of the Minkowski space $\RR^{n-1,1}$.
 
 $\bullet$ {\it Dimension 2.} For $n=2$, the connected  component of the identity in the isometry group is  $SOL \simeq \CC \ltimes \CC^2$, where 
     the action of $\CC$ on $\CC^2$ is given by the complex one-parameter group $I=\left(  \begin{array}{cc}
                                                                 e^t   &   0\\
                                                                 0     &  e^{-t}\\  \end{array}   \right).$  
 
  \subsubsection*{The non-zero constant curvature case}  \label{non-zero curvature}
  The model of the geometry of constant non-zero  curvature,  in dimension $n \geq 2$,
  is the ``holomorphic sphere'' $S_n = O(n+1, {\mathbb C}) / O(n, {\mathbb C})$. Indeed, up a to multiplicative constant, $S_n$  admits a unique,  
  $O(n+1, \CC)$-invariant, holomorphic Riemannian metric $g$. It turns out that $O(n+1, \CC)$  is the full isometry group of $g$,  that $g$  has a constant sectionnal curvature and    is geodesically complete. Therefore,  any $n$-manifold endowed with a holomorphic Riemannian metric of  non-vanishing constant sectionnal curvature  is locally modelled on the geometry $(O(n+1, \CC), S_n)$~\cite{Th}.

 $\bullet$ {\it Dimension 2.}
 A model of $S_2$  is $P^1(\CC) \times P^1(\CC) \backslash Diag$ endowed with the holomorphic Riemannian metric $ \frac{dz_{1}dz_{2}}{(z_{1}-z_{2})^2}$, given in local affine coordinates. Here the isometry group is $SL(2, \CC)$ acting diagonally.

$\bullet$ {\it Dimension 3.}  The unique case where 
$O(n, \CC)$ is not simple is when  $n=4$ and then,  $O(4, \CC) = SL(2, \CC) \times SL(2, \CC)$. 
 The space $S_3$ is identified with the group $SL(2, \CC)$ endowed with a left invariant holomorphic Riemannian metric which  equals the  Killing form at the identity. But the invariance of the Killing form by the adjoint representation implies that this holomorphic Riemannian  metric is also right invariant. Therefore, the right and left multiplicative action of   $SL(2, \CC) \times SL(2, \CC)$ on $SL(2, \CC)$ is isometric.  For more details about this geometry (geodesics...) one can see~\cite{Ghys2}

 \subsection*{Homogeneous spaces}
 
 Left invariant holomorphic Riemannian metrics on a complex Lie group $G$ go down on any compact quotient $\Gamma \backslash G$ by a lattice $\Gamma$. Conversely we have
 the following (see Proposition 3.3 in~\cite{Dum}):
 
 \begin{proposition} Let $g$ be a  holomorphic Riemannian metric on a compact homogeneous space $\Gamma \backslash G$, where $\Gamma$ is a closed subgroup of the complex
 Lie group $G$. Then $\Gamma$ is a lattice  in $G$ and the pull-back of $g$ on $G$ is left invariant.
 \end{proposition}
 
 Note that any  $3$-dimensional unimodular complex Lie group is locally isomorphic to one of the following Lie groups: $\CC^3$, the complex Heisenberg group, the complex $SOL$ group  and $SL(2, \CC)$~\cite{Kir}. 
 
 $\bullet$ $G=\CC^3$. Any left invariant holomorphic Riemannian metric on $\CC^3$ is flat. 
 
 $\bullet$ $G=SL(2, \CC)$. We have seen previously that $SL(2, \CC)$ admits left invariant holomorphic Riemannian metrics of non-zero constant sectionnal curvature.
 
 $\bullet$ $G=Heis$ or $G=SOL$. These groups admit flat left invariant holomorphic Riemannian metrics~\cite{Ra}.

  \subsection*{Nonstandard examples of dimension $3$}

  As above, for  any co-compact lattice $\Gamma$ in $SL(2, \CC)$, the quotient 
  $M=  \Gamma \backslash SL(2, \CC)$  admits a holomorphic Riemannian metric  of non-zero constant sectionnal curvature.
  It is convenient to consider   $M$ as  a quotient of $S_3$ by 
  $\Gamma$,  seen as a subgroup of $O(4, \CC) $ by the 
  trivial embedding $\gamma \in \Gamma \mapsto ( \gamma,1) \in SL(2, \CC) \times SL(2, \CC)$.

New interesting examples of manifolds admitting holomorphic Riemannian metrics of non-zero constant sectionnal curvature have been constructed in \cite{Ghys2}  by deformation of this embedding of $\Gamma$. There, Ghys was interested in  the deformation  of the complex structure of $\Gamma \backslash SL(2, \CC)$,  rather than in  their holomorphic Riemannian metrics. However, one important achievement is the coincidence of complex classification and the holomorphic Riemannian one. 

Examples of deformations of $\Gamma$ 
 are constructed   by means of a morphism 
 $u: \Gamma \to SL(2, \CC)$ and considering the embedding 
 $ \gamma \mapsto ( \gamma, u(\gamma))  $.
 Algebraically, the so obtained action is given 
 by:
 $$(m,\gamma) \in SL(2,\CC) \times \Gamma \to \gamma m u(\gamma^{-1}) \in SL(2,\CC).$$

 It is proved in \cite{Ghys2} that, for $u$ close enough to the trivial morphism,  $\Gamma$ acts properly (and freely)  on $S_3 (\cong SL(2, \CC))$ such that the quotient $M_u$
  is  a complex compact manifold (covered by $SL(2, \CC)$) admitting a holomorphic Riemannian metric of non-zero constant sectionnal curvature. All the $M_u$'s are differentiably diffeomorphic, but are holomorphically diffeomorphic,  iff,  they are isometric (iff, their defining morphisms are conjugate). Note that  left-invariant holomorphic
  Riemannian metrics on $SL(2, \CC)$ which are not right-invariant, in general, will not go down on $M_{u}$.

Let us notice   that despite this systematic study in 
\cite{Ghys2}, there are  still many open questions regarding these examples (including the question of completeness).  A real version of this study is in~\cite{Kulkarni-Raymond, Goldman, Salein}. This story is  also related to the study of Anosov flows with  smooth distributions~\cite{Ghys3}.

 \subsection*{Non-zero constant curvature in higher dimension?}
 
  One interesting 
 problem in differential geometry is to decide if a given homogeneous space $G/I$ possesses or not a compact quotient. A more general related question is to decide
 if there exist compact manifolds locally modelled on $(G,G/I)$ (see, for instance,~\cite{Benoist,  Benoist-Labourie, Labourie, Kobayashi}).

 The case $I = 1$, or more generally $I$ compact,  reduces to the classical question of existence 
 of co-compact lattices in Lie groups.  For  homogeneous  spaces  of non-Riemannian  type (i.e. $I$ non-compact) the problem is much harder. 
 
The  case  
 $S_n = O(n+1, \CC)/ O(n, \CC)$ is a geometric situation where these questions can be tested. It turns out that compact quotients of  $S_n$ are known to exist only for $n=1,3$ or $7$. We discussed the case $n= 3$ above, and the existence of a compact quotient of $S_7$
 was proved  in 
 \cite{Kobayashi}.  Here, we dare ask with \cite{Kobayashi}:
 
 \begin{conjecture} $S_n$ has no  compact quotients,  for $n \neq 1, 3, 7$.

 \end{conjecture}

 A stronger version of our question was  proved in \cite{Benoist}  for $S_n$, if $n$ has the form 
  $4m+1$.
  
  Keeping in mind our geometric approach, we generalize the question 
  to manifolds locally modelled on $S_n$. More exactly:

  \begin{conjecture} A compact complex manifold endowed with a holomorphic Riemannian metric
  of constant non-vanishing curvature is complete. 
   In particular,  such a manifold has dimension 3 or 7. 
  \end{conjecture}

 \section{Geometry of the Killing algebra}

Recall that a holomorphic Riemannian metric $g$ on $M$ is said {\it locally homogeneous} if for all $m,n \in M$ there is a local biholomorphism from an open neighborhood of
$m$ to an open neighborhood of $n$ which sends $m$ to $n$ and preserves $g$. Such a local biholomorphism preserving $g$ is called a {\it local isometry}.

By Theorem~\ref{Dumitrescu}, each holomorphic Riemannian metric on a {\it compact} complex $3$-manifold  is locally homogeneous. Equivalently  the local algebra
of holomorphic Killing fields (i.e. holomorphic vector fields whose local flow preserves $g$) is transitive on $M$. In particular,   the  Killing Lie algebra $\mathcal G$  of $g$ is of dimension $\geq 3$.

Moreover, for any  holomorphic tensor field $\phi$ on $M$, the pseudo-group of local isometries of $g$ preserving also $\phi$ acts transitively on $M$ (i.e.  if we put together $g$ and $\phi$, this  yield to a locally homogeneous geometric structure).

The set of local isometries $I$ of $g$ fixing a point $x_{0} \in M$ generate a local group called  the {\it isotropy group} of $g$. The corresponding Lie algebra $\mathcal I$  consists in the subalgebra of Killing fields vanishing at $x_{0}$. As an isometry fixing $x_{0}$ is uniquely determined by its differential at $x_{0}$~\cite{Wo}, the local group of isotropy at $x_{0}$
injects into the orthogonal group of $(T_{x_{0}}M,g_{x_{0}})$ and thus  it is  of  dimension  $\leq 3$. It follows that 
$\mathcal G$ is of dimension $\leq 6$.

Let $G$ be the connected simply connected complex Lie group corresponding to $\mathcal G$ and $I$ its subgroup corresponding to $\mathcal I$. By a Theorem of Mostow~\cite{Mos},  $I$ is  closed 
in $G$ (this will follow also from our classification of $\mathcal G$ and $\mathcal I$). Thus $g$ is  locally isometric to an algebraic   model $G/I$ endowed with a $G$-invariant holomorphic Riemannian metric.  Since the (full)  isometry group of $G/I$ has at most finitely many connected components, up to a finite cover, $M$ admits a $(G,G/I)$-geometry in Thurston's sense~\cite{Th}: $M$ admits an atlas with open sets in $G/I$ and transition functions given by elements in $G$.

We will classify all possible models  $(G,G/I)$.
We settle first the easiest cases  where $\mathcal G$ has dimension $3, 5$ and $6$.

\subsection{dim $ \mathcal G=3$.}

With Proposition~\ref{unimodulaire} we can easily prove some simplified versions of Theorem \ref{principal}.

     \begin{lemma} \label{2 vector fields}
     
     Let $M$ be a compact connected complex $3$-manifold admitting a holomorphic Riemannian metric $g$.  Assume one of   the following assumptions holds:
     
     (i)   the Killing Lie algebra $\mathcal G$ of $g$  has  dimension $3$;
     
     (ii)  $M$ admits two linearily independent global holomorphic vector fields.
     
    Then,  up to a finite unramified cover, $M$ is a quotient of a complex Lie group $G$ by a lattice $\Gamma$ (hence it  admits some holomorphic Riemannian metric of constant sectionnal curvature) and the pull-back of $g$ on the universal cover of $M$ is  a left invariant holomorphic Riemannian
    metric on $G$.
     \end{lemma}
     
     \begin{remark} If $G=\CC^3$, then $g$ is flat and its  Killing Lie algebra is of dimension $6$ (see Proposition~\ref{dim 6}).
     \end{remark}

     \begin{proof}
     
     (i) As $g$ is locally homogeneous and $\mathcal G$ is of dimension $3$, 
     the action of  $\mathcal G$ on $M$ is simple and transitive. This gives a 
     $(G,G)$-structure on $M$,  where the complex Lie group  $G$  acts on itself by left translations.  The compactness of $M$ implies the completeness of the $(G,G)$-structure~\cite{Th} and  hence 
     $M$ is a quotient of $G$ by a lattice  $\Gamma$.     
     
    (ii) We apply  Theorem \ref{Dumitrescu} to the holomorphic geometric structure on $M$ which is the combination of $g$ with  the two global vector fields . Consequently this
    geometric structure is  locally homogeneous. Moreover,  its Killing Lie algebra is easily seen to be of dimension $3$. Indeed, the local isotropy group at  $x_{0} \in M$  is trivial because
 any element of it which fixes two linearily independent vectors in 
 $T_{x_0}M$ is trivial.  One has just to check directly the claim for the equivalent  situation:        
     $O(3, \CC)$ acting linearily on $\CC^3$.  Finally, we  conclude as in the case (i).
      \end{proof}
      
  \subsection{dim $\mathcal G=6$.}
   
   Here we have the following well-known 
   
   \begin{proposition}    \label{dim 6}
   
   The dimension of $\mathcal G$ is $6$ if and only if $g$ is of constant sectionnal curvature.
   \end{proposition}
   
   \begin{remark}  In this case $\mathcal G$ has a non trivial semi-simple part.
   \end{remark}
   
   \begin{proof}
   
   The dimension of $\mathcal G$ is $6$ if and only if the dimension of $\mathcal I$ is $3$ and if and only if each element in  the   connected component  of identity  of the orthogonal group
   of $(T_{x_{0}}M,g_{x_{0}})$  extends to a local isometry. As the identity component   of the orthogonal group of $(T_{x_{0}}M,g_{x_{0}})$ acts transitively on the set of
   non-degenerate planes in  $T_{x_{0}} M$, all these  planes     have  the same   sectionnal curvature. By local homogeneity,  this sectionnal curvature
   does not  depend on the point $x_{0}$.
   
   Conversely the two models of 3-dimensional spaces of constant sectionnal curvature have a Killing Lie algebra of dimension $6$ which is the Lie algebra of   $O(3, \CC)  \ltimes \CC^3$,
   in the flat case, or the Lie algebra of  $SL(2, \CC) \times SL(2, \CC)$, in the non flat one.
   \end{proof}
   
\subsection{dim $\mathcal G=5$.}
We will see this   never happens.

Recall first  that $SL(2, \CC)$ is locally isomorphic to $O(3, \CC)$. One way to see it is to consider the adjoint representation of $SL(2, \CC)$ into the 3-dimensional complex vector space $sl(2, \CC)$ and to note that this action preserves the Killing form. More precisely,  we have $SO(3, \CC) \simeq PSL(2, \CC)$, where 
$SO(3, \CC)$ is the connected component of the identity of the orthogonal group and $PSL(2, \CC)$ is the quotient of $SL(2,\CC)$ by  the center $\{ Id, -Id \}$.

\begin{proposition}
The dimension of $\mathcal G$ is $\neq 5$.
\end{proposition}

\begin{proof}

Assume, by contradiction,  that   dim $\mathcal G = 5$ and,  equivalently,  the dimension of the isotropy  $\mathcal I$ is $2$. Consider the action of the local isotropy group at $x_{0}$
on $T_{x_{0}}M$ and identify this local isotropy to a 2-dimensional subgroup $I$ of $SO(3, \CC) \simeq  PSL(2, \CC)$. The action of $I$ on $T_{x_{0}}M$ preserves $g_{x_{0}}$, but also the curvature tensor and, in particular,
the Ricci tensor $Ricci_{x_{0}}$ which is a complex quadratic form on $T_{x_{0}}M$.

Consider the action of $PSL(2, \CC)$  on the complex vector space of complex quadratic forms $S^2(T^*_{x_{0}}M)$.
This action preserves $g_{x_{0}}$ and gives  an action of  $PSL(2, \CC)$ on the quotient vector space  $S^2(T^*_{x_{0}}M) /  \CC g_{x_{0}}$. 
  
  The isotropy group lies in the stabilizer of the class of  $Ricci_{x_{0}}$ in the quotient $S^2(T^*_{x_{0}}M) /  \CC  g_{x_{0}}$. But, for an algebraic  action of $PSL(2, \CC)$
  on an affine space,  the stabilizer of an element can not be 
    1-dimensional. Indeed, by contradiction,  up to an inner automorphism of $PSL(2, \CC)$, the stabilizer  coincides with   the subgroup $G' \subset PSL(2, \CC)$ of upper triangular matrices and  thus the orbit $PSL(2, \CC) /G'$ is biholomorphic to the projective line $P^1(\CC)$,  which is compact and so can not be holomorphically embedded in an affine space.
  
  It follows that the stabilizer of the $Ricci_{x_{0}}$ class in $S^2(T^*_{x_{0}}M) /  \CC  g_{x_{0}}$ is of dimension $3$ and hence equal to $PSL(2, \CC)$. This implies that   
   $Ricci_{x_{0}} = \lambda g_{x_{0}}$, with  $\lambda \in \CC$   and the function    $\lambda$ is constant on $M$ by local homogeneity. But then, $g$ has constant
   sectionnal curvature and so $\mathcal G$ is of dimension $6$ which is contrary to our initial assumption.    
 
   \end{proof}

   \subsection{dim $\mathcal G= 4$ } This is the most delicate case and all our analysis throughout the paper will devoted to it.

   Here  $\mathcal I$ has dimension $1$.  The (local) isotropy group $I$ is algebraic and has finitely many components. Up to a finite cover, we can assume it connected, i.e. a one parameter group. Therefore, $ I$  is conjugate in $PSL(2, \CC)$ to 
  one of the following:
  
  \begin{enumerate}

\item   A  {\it unipotent}  one-parameter subgroup $\left(  \begin{array}{cc}
                                                                 1   &   t \\
                                                                 0     &  1\\
                                                                 \end{array} \right)$  fixing in $T_{x_0}M$ a vector of norm   $0$;

\item   A  {\it semi-simple} one-parameter subgroup    $\left(  \begin{array}{cc}
                                                                 t   &   0\\
                                                                 0     &  t^{-1} \\
                                                                 \end{array} \right)$ fixing in $T_{x_0}M$  a vector of norm   $1$.

\end{enumerate}

\subsection*{Adapted basis}
   
   In order to understand  the action of $I$ on $T_{x_{0}}M$ (as a subgroup of $O(3, \CC)$) we shall consider some adapted bases.
   
   Let us  first consider the case where    the isotropy is semi-simple. Then the action of $I$ on $T_{x_{0}}M$ fixes some vector  $e_{1}$ of norm   $1$. The plane  $e_{1}^{\bot}$ is non degenerate
   and, up to a multiplicative constant, the vectors  $e_{2}, e_{3} \in e_{1}^{\bot}$ are uniquely defined by the following conditions: $e_{2}, e_{3}$ generate the two isotropic directions in $e_{1}^{\bot}$
   and $g(e_{2}, e_{3})=1$. The time $t$ of the flow generated by the isotropy $\mathcal I$ will be given in this adapted basis $(e_{1},e_{2}, e_{3})$, by the formula $(e_{1}, e_{2}, e_{3}) \to (e_{1}, e^t e_{2}, e^{-t} e_{3})$.
   
   In the case of a unipotent isotropy,  the action of  $I$ on  $T_{x_{0}}M$ fixes an isotropic vector $e_{1}$ and so preserves the degenerate plane $e_{1}^{\bot}$ (of course $e_{1} \in e_{1}^{\bot}$). 
  In order to define an adapted basis,   take two vectors $e_{2}, e_{3} \in T_{x_{0}}M$ such that:
    $g(e_{1}, e_{2})=0$, $g(e_{2},e_{2})=1$, $g(e_{3},e_{3})=0,
 g(e_{2},e_{3})=0$ and  $g(e_{3}, e_{1})=1.$ 
 
  Note that such an adapted basis is uniquely determined by the choice of an unitary vector $e_{2} \in e_{1}^{\bot}$. Indeed, then $e_{3}$ is   uniquely defined in 
 $e_{2}^{\bot}$ by the relation  $g(e_{3}, e_{1})=1$ ( $e_{1}$ and  $e_{3}$ generate the two isotropic directions in $e_{2}^{\bot}$). 
 
 The action of the isotropy $ I$ on $T_{x_{0}}M$ sends an adapted basis to  an adapted basis. This action is given in the basis $(e_{1}, e_{2}, e_{3})$ by 
 $\left(  \begin{array}{ccc}
                                                                 1   &   t & - \frac{t^2}{2}\\
                                                                 0     &  1 &  -t\\
                                                                 0     &   0  &  1\\
                                                                 \end{array}  \right).$

   \begin{lemma}  \label{non trivial center}
   
 (i)  If $\mathcal G$ is of dimension $4$, then, up to a finite cover,  $M$ admits a global holomorphic vector field $X$ which is preserved by the action of $\mathcal G$. The norm of  $X$
   is constant equal to $0$ or constant equal to $1$, according to that  the isotropy is unipotent or semi-simple.
   
 (ii) The divergence of $X$ (with respect of  the volume form of $g$)  is $0$.
   
 (iii) If the isotropy is semi-simple,  then $X$ is a  Killing field.
 \end{lemma}
 
 \begin{corollaire}
 
 If the isotropy is   semi-simple, then $\mathcal G$ has a non trivial center.
 \end{corollaire}
 
 \begin{proof}
 
 (i) At $x_0$, $X$ is defined by $X(x_0) = e_1$.

 (ii)  Denote by $\phi^t$ the complex flow generated by $X$. Recall that the divergence $div(X)$ of $X$, with  respect to the volume form $vol$ of $g$, is given by the formula
 $L_{X}vol=div(X)vol$, where $L_{X}$ is the Lie derivative in the direction $X$. As $\mathcal G$ acts transitively on $M$ preserving $X$ (and also $vol$),  the function $div(X)$
 is holomorphic and so is a  constant $\lambda \in \CC$. This means that $(\phi^t)^*vol=e^{\lambda t}vol$, for all $t \in \CC$. But the total real volume of $M$ given by the integral on $M$
 of the real form $vol \wedge \overline{vol}$ has to be preserved by  $\phi^t$. Thus the modulus of $e^{\lambda t}$ equals  $1$ for all $t \in \CC$. It then follows that $\lambda=0$, that is
  $div(X) = 0$. 
 
 (iii) The action of $\mathcal G$ preserves $X$ and so also $X^{\bot}$. We will show first that $\phi^t$ preserves $X^{\bot}$ as well.
 Take a point $x_{0} \in M$ and consider its image  $\phi^t(x_{0})$. For each  $t \in \CC$ let us choose 
    a local isometry  $g^t$ sending  $x_{0}$ to $\phi^t(x_{0})$.
    
    The local diffeomorphism $(g^t)^{-1} \circ \phi^t$ fixes  $x_{0}$ and the vector  $X(x_{0}) \in T_{x_{0}}M$. Since $X$ is $\mathcal G$-invariant, $(g^t)^{-1} \circ \phi^t$ commutes with all local isometries . In particular,  the differential $L_{t}$ of  $(g^t)^{-1} \circ \phi^t$
    at  $x_{0}$ commutes  with the action of the isotropy at  $x_{0}$ and hence preserves   the eigenspaces of the isotropy.
   Since the isotropy is supposed to be semi-simple,   the differential $L_{t}$ preserves  the non-degenerate plane $X(x_{0})^{\bot}$ and also its two isotropic directions. 
   
   As $div(X)=0$,  the differential  $L_{t}$ preserves the volume. It follows that the product   of the two eigenvalues corresponding to the two isotropic directions of $X(x_{0})^{\bot}$ equals $1$. This implies that the differential of $(g^t)^{-1} \circ \phi^t$ at  $x_{0}$ is an isometry. Consequently the flow of $X$ acts by isometries and $X$ is Killing. Hence $\CC X$ is in the center of $\mathcal G$.
   \end{proof}

     \begin{proposition}   \label{invariance}

  If the isotropy is unipotent,  then the holomorphic field of  complex endomorphisms $\nabla_{\cdot }     X$ of  $TM$,    in an adapted basis, is  $ \left(  \begin{array}{ccc}
                                                                 0    &   0 & \alpha \\
                                                                 0     &  0  &  0  \\
                                                                 0     &   0  &  0 \\
                                                                 \end{array} \right), $ with $\alpha$ a complex constant.
                                                                 
Then $X$ is Killing if and only if $\alpha=0.$
                                                                 
      \end{proposition}

 \begin{proof}
 Let us fix   $x_{0} \in M$ and an adapted basis $(e_{1}, e_{2}, e_{3})$ of  $T_{x_{0}}M$.       In this basis the differential $L_{t}$ of $I$ at $x_{0}$   is given by the one-parameter group
 $\left(  \begin{array}{ccc}
                                                                 1   &   t & - \frac{t^2}{2}\\
                                                                 0     &  1 &  -t\\
                                                                 0     &   0  &  1\\
                                                                 \end{array}  \right). $

       First we show that  any  $\mathcal G$-invariant holomorphic field of  complex endomorphisms $\Psi$ of $TM$ has,  in our adapted basis,  the following  form:                                                               
$  \left(  \begin{array}{ccc}
                                                                 \lambda    &   \beta  & \alpha \\
                                                                 0     &  \lambda  &  - \beta  \\
                                                                 0     &   0  &  \lambda \\
                                                                 \end{array} \right),$      with $\alpha, \beta$ and $ \gamma \in \CC.$
                                                                 
      Let  $B$ be the matrix of $\Psi(x_{0})$ in the basis  $(e_{1}, e_{2}, e_{3})$. Since $\Psi$ is $\mathcal I$-invariant,  $B$ and  $L_{t}$ commute. 
      Each eigenspace of $B$ is  preserved by $L_{t}$ and conversely. As $L_{t}$ does not  preserve any non trivial spliting 
      of $T_{x_{0}}M$, it follows that all eigenvalues of $B$ are equal to some $\lambda \in \CC$.    A straightforward calculation shows that $B$ has  the previous form.
      As $\Psi$ is $\mathcal G$-invariant,  the parameters  $\alpha, \beta$ and $\gamma$    do not 
      depend of $x_{0}$.
      
      We apply this result to $\nabla_{\cdot} X$     (which is $\mathcal G$-invariant because $X$ and $\nabla$ are). As the trace of $\nabla_{\cdot} X$   is the divergence of $X$,
       lemma~\ref{non trivial center} implies  $\lambda =0$. 
      
      It will be (independently) shown  in Proposition \ref{X parallel}  that $X$ is parallel on any direction  tangent to $X^{\bot}$. It  follows that
      $\nabla_{e_{2}} X  =0$ and  $\beta=0$.
      
      The vector field $X$ is Killing if and only if $\nabla_{\cdot} X$ is $g$-skew-symmetric~\cite{Wo}. But an endomorphism of rank $\leq 1$ is skew-symmetric if and only if
      it is trivial. It follows that $X$ is Killing if and only if $\alpha=0$.      
    \end{proof}

  \subsection*{Geodesic  foliations}
   The following lemma is just the complexification in the realm of holomorphic Riemannian metrics of a well-known fact  remarked for the first time by M. Gromov \cite{Gro}
 (see also the survey \cite{DG})  in the context of  Lorentz geometry.
  
  \begin{lemma} 
  
  (i) If the isotropy is unipotent,  then the plane field $X^{\bot}$ is integrable. Its  tangent holomorphic foliation of codimension one $\mathcal F$   is geodesic, $g$-degenerate  and $\mathcal G$-invariant.
  
  (ii) If the isotropy is semi-simple,  then $M$ possesses two holomorphic foliations of codimension one $\mathcal F_{1}$ and $\mathcal F_{2}$, which are geodesic, $g$-degenerate
  and $\mathcal G$-invariant. The tangent space of each one of these two foliations is generated by $X$ and by one of the two isotropic directions of $X^{\bot}$.
  \end{lemma}
  
  \begin{proof}
  
  The idea of Gromov's proof is to consider the graph of a local  isometry fixing $x_{0} \in M$ as a (3-dimensional) submanifold in $M \times M$ passing through $(x_{0},x_{0})$.
  This submanifold is geodesic and isotropic for the holomorphic Riemannian metric $g \oplus (-g)$ on $M \times M$. If $f_{n}$ is a sequence of elements in the local  isotropy
  group at $x_{0}$ (identified with the orthogonal  group of $(T_{x_{0}}M, g_{x_{0}})$) which goes to  infinity in this orthogonal group, then the sequence  of corresponding
  graphs tends to a 3-dimensional geodesic and isotropic  submanifold $F'$ which  is no longer a graph. Nevertheless,   the intersection of $F'$ with the vertical space 
  $\{ x_{0}  \}  \times M$ is isotropic in $M$ and thus has dimension $\leq 1$. The projection $F$ of $F'$ on the horizontal space $M \times \{ x_{0}  \}$ is a geodesic surface
  passing through $x_{0}$.
  
  In our situation $I$ has dimension $1$ and we can take a sequence of elements of the one-parameter group $I$ in  the orthogonal group going to infinity (one parameter groups are not  compact, which contrasts with the real case). In exponential 
  coordinates our local isometries are linear and in some adapted basis they have the form presented previously. We note that the limit of our sequence of (linear) graphs 
  is the plane $X(x_{0})^{\bot}$ if the isotropy is unipotent and the two planes generated by $X(x_{0})$ and by each of the two isotropic directions of $X(x_{0})^{\bot}$ if 
  the isotropy is semi-simple.
  
  These foliations are obviously $\mathcal G$-invariants, as everything is. 
   \end{proof}

  We will also denote by $X$ and $\mathcal F$ the corresponding vector field and foliation on the algebraic model $G/I$.

  \subsection*{The stabilizer $H$ of a leaf}
  
  If the isotropy is unipotent, denote by $\mathcal H$ the subalgebra of $\mathcal G$ stabilizing the leaf  $\mathcal F(x_{0})$ of $\mathcal F$ passsing through $x_{0} \in M$ and by $H$ the corresponding Lie subgroup of $G$. We keep the same
  notation for the stabilizer of $\mathcal F_{1}(x_{0})$ if the isotropy is semi-simple.
  
  \begin{proposition}  \label{stabilizer}
  The group $H$ is of dimension $3$ and acts transitively on $\mathcal F(x_{0})$ (or $\mathcal F_{1}(x_{0})$ accordingly). The isotropy $I$ at  $x_{0}$ lies in $H$.
  \end{proposition}
  
  \begin{corollaire} The leaf $F$ is locally modelled on $(H,H/I)$.
  \end{corollaire}
  
  \begin{proof}
  
  We give the proof in the case of unipotent isotropy.
   Take $x_{1} \in \mathcal F(x_{0})$ and consider a local isometry $\phi$ sending $x_{0}$ on $x_{1}$. As $\phi$  preserves $X$ and $X^{\bot}$ it has to send $exp_{x_{0}}(X^{\bot})$
  onto  $exp_{x_{1}}(X^{\bot})$. The leaf $\mathcal F(x_{0})$ being geodesic, $exp_{x_{0}}(X^{\bot}) \subset \mathcal F(x_{0})$ and $exp_{x_{1}}(X^{\bot})  \subset \mathcal F(x_{0})$.
  That means that $\phi$ lies in the stabilizer of $\mathcal F(x_{0})$. In particular,  if $\phi$ fixes $x_{0}$ then $\phi$ lies  in the stabilizer of $\mathcal F(x_{0})$. This implies
  $\mathcal I \subset \mathcal H$.
  
  As $\mathcal G$ acts transitively on   $\mathcal F(x_{0})$, the previous argument shows that $\mathcal H$ acts transitively on   $\mathcal F(x_{0})$   (with isotropy of dimension $1$).
  It follows that $\mathcal H$ has dimension $3$.
  \end{proof}

\section{Algebraic models for the local structure: the semi-simple case}   \label{algebraic models semi-simple}

In this section  the Killing algebra $\mathcal G$  has  dimension $4$, and  thus the isotropy $\mathcal I$ has dimension $1$.
 We assume that $\mathcal G$ has a non-trivial semi-simple part.  
 
 \begin{proposition}

 Assume $\mathcal G$ has a non-trivial  semi-simple part. Then, it is a direct product of Lie algebras $\CC \oplus sl(2, \CC)$, and we have two possible models $G/I$:

(1) The holomorphic Riemannian metric is  left invariant  on the group   $SL(2, \CC)$.

 The identity connected component of  its  isometry group   is a direct product of  $SL(2,\CC)$  acting 
by left translations  and  some one parameter subgroup $h^t \subset SL(2, \CC)$ acting on by right translations.  The isotropy group $I$ is the image of the diagonal embedding $(h^t,h^t)$  in $\CC \times SL(2,\CC)$.

(2) The holomorphic Riemannian direct product $\CC \times S_{2}$,  where $S_{2}$ is the universal  model  of  a  surface with holomorphic Riemannian metric of  non zero constant sectionnal curvature and 
  $\CC$ is endowed with its  standard  metric  $dz^2$.
  
 The action of  the isometry group  $G=\CC \times SL(2,\CC)$ is split. 
The isotropy $I$  is the  one-parameter subgroup of $SL(2,\CC)$ given by
 $\left(  \begin{array}{cc}
                                                                 e^t   &   0\\
                                                                 0     &  e^{-t}\\
                                                                 
                                                                 \end{array}  \right).$

 \end{proposition}

\begin{corollaire}
In the case $(1)$  the action of $\mathcal G$ on $M$ preserves the holomorphic Riemannian metric of non zero constant sectionnal curvature coming from the Killing form on $sl(2, \CC)$.
\end{corollaire}

\begin{remark} It will be shown in \S\ref{semi-simple isotropy} that 
the situation $(2)$ cannot occur  on  compact $3$-manifolds.
\end{remark}

\begin{proof}

There is no semi-simple algebra of dimension 4, and $sl(2, \CC)$ is the unique semi-simple complex Lie algebra of dimension 3. Therefore, ${\mathcal G}$ is a direct 
product $\CC \times sl(2, \CC)$ (see, for instance,~ \cite{Kir}). 

If the isotropy of some point  intersects non-trivially the factor $SL(2, \CC)$, then this is the case for all points. In fact, since the isotrop $I$ has dimension 1, it intersects 
$SL(2, \CC)$ iff it is contained inside it.  \\

(1) Therefore, in the case of trivial intersection, the group 
$SL(2, \CC)$ acts freely transitively on $M$. The metric is thus identified to a left invariant one on $SL(2, \CC)$. 

Consider the action of the isotropy $I$ on $SL(2, \CC)$ (the base point being the neutral element $Id$  in $SL(2, \CC)$). Our claim  reduces to the fact that the $I$-action  coincides with the adjoint action of some one parameter group $h^t$. 
For this, it suffices to show that the  metric is preserved by the adjoint action of $h^t$ on 
$sl(2, \CC)$. Indeed, if so,  this integrates on the adjoint action of $h^t$
on the group $SL(2, \CC)$ which is isometric. But, since the dimension of the isotropy  is one, we get coincidence of $I$ with  the adjoint action of $h^t$.

The $I$-action  on $sl(2,\CC)$ by the adjoint representation is done
    by Lie algebras isomorphisms.
    
    On the other hand the previous action identifies with the $I$-action on $T_{Id}SL(2,\CC)$ and  has to  fixe  some vector. It is easy to check that each one-parameter group of isomorphisms of the Lie algebra $sl(2, \CC)$ fixing a vector coincides with  the adjoint representation of some one-parameter subgroup $h^t$ of $SL(2, \CC)$.

 (2) Assume now that   $ I \subset SL(2, \CC)$.  The action of $I$ on $\CC \oplus sl(2, \CC)$ gives an $I$-invariant non trivial splitting of $T_{x_{0}}M$. It follows that $I$
 is semi-simple and the $SL(2, \CC)$-orbits are tangents to $X^{\bot}$ (in particular, they are $g$-non degenerate).
 Then, the $SL(2, \CC)$-orbits are complex homogeneous surfaces endowed with a $SL(2, \CC)$-invariant holomorphic Riemannian metric. They have in particular constant curvature, and obviously cannot be flat (because their  Killing algebra contains $sl(2, \CC)$). Up to a multiplicative constant, they are isometric to $S_{2}$.
\end{proof}

\section{Algebraic models for the local structure: the solvable  case}   \label{algebraic models}

We assume here that $G$ is solvable (and of  dimension 4). 

\begin{proposition}    \label{stabilizer} (i) The derivative Lie algebra $\lbrack \mathcal{H}, \mathcal{H}  \rbrack$ is 1-dimensional.

(ii)The group $H$ is isomorphic either to the Heisenberg group or to the product $\CC \times AG$, where $AG$ is the universal covering of the affine group  of the complex line.
\end{proposition}

Recall  that the affine group of the complex line  is the group of transformations of $\CC$, given by $z \to az+b$, with $a \in \CC^*$ and $b \in \CC$. If $Y$ is the infinitesimal generator of the homotheties and $Z$ the infinitesimal generator of the translations,  then $\lbrack Y, Z \rbrack =Z$. 
\begin{proof}

(i) It is a general fact that a derivative algebra of a solvable algebra is nilpotent. Remark first that $\lbrack \mathcal H, \mathcal H \rbrack \neq 0$. Indeed,   if not $\mathcal H$ is abelian and the action of the isotropy $\mathcal I \subset \mathcal H$ would be trivial
on $\mathcal H$ and hence on $T_{x_{0}}F$ which is identified to $\mathcal H / \mathcal I$. Since the restriction to the isotropy action to the tangent space of $F$ is injective
this  implies that the isotropy action is trivial on $T_{x_{0}}G/I$ which is impossible.

As $\mathcal H$ is 3-dimensional, its derivative algebra $\lbrack \mathcal H, \mathcal H \rbrack$
is a nilpotent Lie algebra of dimension $1$ or $2$, hence $\lbrack \mathcal H, \mathcal H \rbrack \simeq \CC$ or $\lbrack \mathcal H, \mathcal H \rbrack \simeq \CC^2$.

Assume by contradiction that $\lbrack \mathcal H, \mathcal 
H \rbrack \simeq \CC^2$. \\

  We first  prove  that the isotropy $\mathcal I$ lies in $\lbrack \mathcal H, \mathcal H \rbrack$. If not, $\lbrack \mathcal H, \mathcal H \rbrack \simeq \CC^2$ will act freely and so transitively
on $F$. Therefore $F$ is identified with the group $\CC^2$ endowed with a left invariant connection, a left invariant holomorphic degenerate Riemannian metric (compatible with the connection) and a left invariant 
holomorphic vector field (which is $X$). 

We  show now that the connection is flat.
The local model for  the left invariant degenerate metric on $F$ 
is  $dh^2$ in the coordinates $(x,h)$ of $\CC^2$. In this coordinates the left invariant vector field $X$ coincides with $\frac{\partial}{\partial x}$,
if the isotropy is unipotent and with $\frac{\partial}{\partial h}$, if the isotropy is semi-simple. 

An easy calculation shows  that any torsion-free and $\CC^2$-invariant connection compatible with $dh^2$ is given by $\nabla_{\frac{\partial}{\partial h}}\frac{\partial}{\partial h}=a
\frac{\partial}{\partial x}$, $\nabla_{\frac{\partial}{\partial x}}\frac{\partial}{\partial x}=b
\frac{\partial}{\partial x}$ and $\nabla_{\frac{\partial}{\partial h}}\frac{\partial}{\partial x}=
\nabla_{\frac{\partial}{\partial x}}\frac{\partial}{\partial h}=c
\frac{\partial}{\partial x}$, for some   $a,b,c \in \CC$.
The invariance by the isotropy one-parameter group implies that at least two of the parameters $a,b,c$ vanish. In this case the curvature of $\nabla$ vanishes.

The isometry group of this model is $\CC \ltimes \CC^2$, where the action of the isotropy $I \simeq \CC$ on $\CC^2$ is given by the one parameter  group of linear transformations
$\left(  \begin{array}{cc}
                                                                 1   &   t\\
                                                                 0     &  1\\
                                                                 
                                                                 \end{array}  \right) $,  if $I$ is unipotent, or by $\left(  \begin{array}{cc}
                                                                 e^t   &   0\\
                                                                 0     &  1\\
                                                                 
                                                                 \end{array}  \right) $, 
                                                                if $I$ is semi-simple. Our group is thus   
                                                                isomorphic to                                        
    the  Heisenberg group or to  $AG \times \CC$. In both cases the derivative group is 1-dimensional which  contradicts  our assumption, and hence  $\mathcal I \subset \lbrack \mathcal H, \mathcal H \rbrack$. \\

            It follows in particular that the orbits of $\lbrack \mathcal H, \mathcal H \rbrack$     on $F$ are 1-dimensional. 
          We prove now that the orbits of      $\lbrack \mathcal H, \mathcal H \rbrack$     on $F$ correspond to the isotropic direction in $F$  and the isotropy $I$ is unipotent.
           
           Let $Y$ be a generator of $\mathcal I$, $\{Y,X'  \}$ be generators of $\lbrack \mathcal H, \mathcal H \rbrack$ and   $\{Y,X' ,Z \}$      be a basis   of $\mathcal H$. The tangent
           space of $F$ at some point $x_{0}  \in F$ is identified with $\mathcal H /  \mathcal I$ and the infinitesimal  (isotropic) action of $Y$ on this tangent space is given      in the basis
           $\{ X', Z \}$ by the matrix $ad(Y) = \left(  \begin{array}{cc}
                                                                 0  &   *\\
                                                                 0     &  0\\
                                                                 
                                                                 \end{array}  \right) $. This is because $\lbrack \mathcal H, \mathcal H \rbrack  \simeq \CC^2$ and $ad(Y)(\mathcal H) \subset \lbrack \mathcal H, \mathcal H \rbrack.$
                        Moreover, $ad(Y) \neq 0$ since the restriction to the isotropy action to the tangent space of $F$ is injective. 
                        
                        From this form of $ad(Y)$,  we see that the isotropy is unipotent with fixed direction $\CC X'$. This direction is exactly the tangent direction of the orbits of 
$\lbrack \mathcal H, \mathcal H \rbrack$ on $F$. \\

Denote by $\mathcal L$ the derivative algebra $\lbrack \mathcal G, \mathcal G \rbrack $ of  $\mathcal G$. 
Then $\mathcal L \supset \lbrack \mathcal H, \mathcal H  \rbrack \supset \mathcal I.$ The dimension of $\mathcal L$ is $2$ or $3$ and the $\mathcal L$-orbits on $G/I$ have dimension $1$ or $2$ accordingly. \\

Assume first that $\mathcal L$ is 3-dimensional and thus has 2-dimensional orbits on $G/I$. The foliation of $G/I$ provided by the $\mathcal L$-action is 2-dimensional
and invariant by the unipotent isotropy $I$. Since  $X'^{\bot}$ is the only plane field on $G/I$ preserved by the isotropy, it follows that the leafs of the $\mathcal L$-action
coincide  with those of the $\mathcal H$-action. So $\mathcal L = \mathcal H$, as Killing algebra of $F$. But this is impossible, since  $\mathcal L$ is nilpotent (as a derivative algebra of a solvable algebra) and ${\mathcal H}$ is not (its derivative algebra is supposed  to be $2$-dimensional).\\

It remains to settle the case where $\mathcal L$ is 2-dimensional. We show in this case that the infinitesimal isometry $ad(Y)$ of $T_{x_0}G/I$ has rank 1, which is 
not possible for an infinitesimal isometry of a holomorphic Riemannian metric.

Since $\mathcal L =\lbrack \mathcal G, \mathcal G \rbrack$, the image of $\mathcal G$ by the isotropy
action $ad(Y)$ at $x_{0} \in G/I$  is contained in $\mathcal L$. Thus this image has at most dimension  $2$ and as $\mathcal I $$ \subset $$ \mathcal L$ and the tangent space at $x_{0}$ is identified 
with $\mathcal G /  \mathcal I$,  the image of $ad(Y)$ in $T_{x_{0}} G/I$ is of dimension at most $1$.

 This completes the proof of   part $(i)$ of the proposition. \\

(ii) Let $Z$ be a generator of $\lbrack \mathcal H, \mathcal H \rbrack$ and consider its adjoint map $ad(Z) : \mathcal H \to \CC Z$. If this map is trivial then,  $Z$ is central
  and $\mathcal H$ is nilpotent isomorphic to the Heisenberg group.

Consider now the case where $ad(Z)$ is not  trivial. Let $X'$ be a generator of the kernel of $ad(Z)$  and take $Y \in \mathcal H$ such that $\{ Y, X', Z \}$ is a basis of
$\mathcal H$. We can assume that $\lbrack Y, Z \rbrack =Z$. We also have $\lbrack X', Y \rbrack =aZ$, with $a \in \CC$. After replacing 
 $X'$ by $X' +aZ$,  we can assume that 
$a=0$. It follows that $H=\CC \times AG$,  where the center of $H$ is  $exp(\CC X')$ and $AG$ is generated by $exp(\CC Z)$ and $exp(\CC Y)$.
\end{proof}

\subsection{The case: $H=\CC \times AG$}

In this case, all possible algebraic  models $(G,G/I)$ are described in the following:

\begin{proposition}  \label{semi-simple classification}
   
The isotropy group $I$ is semi-simple (it is generated by the infinitesimal generator of the homotheties in $AG$) and   $G$ is one of the following Lie groups:

\begin{enumerate}
\item
                                    $G= \CC \times SOL$

                                    \item

                                     $G = \CC \ltimes Heis$ 
                                     
                                     \item 
                                     
                                     $G = \CC^2 \ltimes \CC^2$
      \end{enumerate}

      In  case $(2)$ the action of the first factor $I \simeq \CC$ on $Heis$, 
        is given by $(X',Z,T) \to (X', e^tZ,e^{-t}T)$,    
      with  respect of a basis $(X',Z,T)$, such that  $X'$ is central and   $\lbrack T, Z \rbrack   =X'$,

       In case $(3)$ the action of the first copy  of $\CC^2$ on the second one is given by the matrices     $\left(  \begin{array}{cc}
                                                                 e^t   &   0\\
                                                                 0     &  e^{-t}\\
                                                                 
                                                                 \end{array}  \right) $  and  $\left(  \begin{array}{cc}
                                                                 1   &   0\\
                                                                 0     &  e^{-t}\\
                                                                 
                                                                 \end{array}  \right) $.                                                                 
 \end{proposition}
 
\begin{remark}
  As the center of $\mathcal G = \CC^2 \ltimes \CC^2$  is trivial, it follows from Lemma~\ref{non trivial center} that  this Lie algebra cannot occur as a local Killing algebra
 for a  holomorphic Riemannian metric on a   compact complex 3-manifold.
\end{remark}

\begin{proof}

As before suppose that $\{X',Y,Z \}$ is a basis of $\mathcal H$ with $X'$  central   and $Y,Z$ spanning the Lie algebra of $AG$ such that  $\lbrack Y,Z \rbrack =Z$.
Denote by $T$ a fourth generator of the Killing algebra $\mathcal G$.

We show that, up to an automorphism of $\mathcal H$ sending $Y$ to $Y+aZ+bX'$, with $a,b \in \CC$, {\it the isotropy algebra $\mathcal I$ is $\CC Y$.}

Observe that  $ad(\alpha X' + \beta Z)(\mathcal H) \subset \CC X' \oplus  \CC Z$, for all $\alpha, \beta \in \CC$. If the isotropy $\mathcal I$ is $\CC (\alpha X' + \beta Z)$ then
the action  of $ad(\alpha X' + \beta Z)$ on  $T_{x_{0}}F \simeq \mathcal H /  \mathcal I$ is given by a matrix of rank $1$. Consequently the isotropy is not semi-simple.
We then   proved that in the case where the isotropy is semi-simple, the isotropy $\mathcal I$ doesn't lie in $\CC X' \oplus  \CC Z$ and, up to an automorphism of $\mathcal H$
sending $Y$ to $Y+aZ+bX'$, we can assume that $\mathcal I = \CC Y$.

Now,  we show the same result in the case of  unipotent isotropy. Observe  first that $\mathcal I \neq \CC X'$ since  the central element $X'$ acts trivially on $\mathcal H$
and hence also on $\mathcal H / \mathcal I \simeq T_{x_{0}}F$,  which is impossible.

Assume, by contradiction,  that $\mathcal I \subset \CC X' \oplus \CC Z$. Up to an automorphism of $\mathcal H$ sending $Z$  to $Z + \alpha X'$, with $\alpha \in \CC$, we can assume
that $\mathcal I = \CC Z$. Then,  the abelian Lie algebra $\CC X' \oplus \CC Y$ intersects trivially $\mathcal I$ and  will act freely and transitively on $F$. As in the proof of
Proposition~\ref{stabilizer}, this implies that $F$ is flat and  the                                                              
Killing Lie  algebra of $F$ is $heis$. But,  this is impossible,   since  the Heisenberg group is nilpotent and $H=\CC \times AG$ is not.
                                                                 
  It follows that, up to an automorphism of $H$, we have  $\mathcal I = \CC Y$. This is impossible in the  unipotent isotropy case. Indeed,  the abelian Lie algebra $\CC X'  \oplus \CC Z$
              acts freely and transitively on $F$ and $F$ is flat. If the isotropy was unipotent then, as before,  $H$ is  isomorphic to the Heisenberg group which  contradicts  our hypothesis.
              
             Therefore, the isotropy is semi-simple. As the isotropy $\CC Y$ fixes $X'$ and expands the direction $\CC Z$ (because of the relation $\lbrack Y, Z \rbrack =Z$), we can choose
             as fourth generator $T$ of $\mathcal G$ the second isotropic direction of the Lorentz plane $X'^{\bot}$. Then we will have $\lbrack Y, T \rbrack =-T + aY$, for some
             constant $a \in \CC$ and we can replace $T$ with  $T -aY$ in order to get $\lbrack Y,T \rbrack =-T$.
         
        In the following, we  assume that $\lbrack Y,T \rbrack =-T$.
       
             We will first  show that               
              $\lbrack T,Z \rbrack =aX' + bY$, with  $a,b \in \CC$ and  $\lbrack T, X' \rbrack=cT$, for some $c \in \CC$.
     
     For the first relation we use the Jacobi relation  $\lbrack Y, \lbrack T, Z \rbrack \rbrack=  \lbrack  \lbrack Y, T \rbrack , Z \rbrack +  \lbrack T, \lbrack Y, Z \rbrack \rbrack =
     \lbrack -T, Z \rbrack + \lbrack T, Z \rbrack =0$ to get that $ \lbrack T, Z \rbrack $ commutes with  $Y$ and consequently lies in  $\CC Y \oplus \CC X'$.
     
     To get the second one, observe   that  $X'$ et $Y$ commute, and thus  $T$ (which is an eigenvector of  $ad(Y)$, is also  an eigenvector of $ad(X')$). This gives
     $\lbrack T, X' \rbrack= c T$, for some  $c \in \CC$. \\
          
      Consider now the derivative algebra $\mathcal L = \lbrack \mathcal G, \mathcal G \rbrack$ and recall it is nilpotent.

      The relations    $\lbrack Y, Z \rbrack = Z$,        $\lbrack Y, T \rbrack = -T$ and       $\lbrack T,Z \rbrack =aX' + bY$, 
      show that $\mathcal L$ contains the Lie algebra generated by  $Z,T$ and  $aX'+bY$.      We have $\lbrack aX' + b Y, Z \rbrack = bZ$ and this implies  $b=0$ (if not 
      the Lie algebra generated by  $aX'+ b Y$ and $Z$  is isomorphic to the Lie algebra of $AG$, which is not nilpotent and so  cannot be embedded  into the nilpotent algebra $\mathcal L$).
        It follows that $b=0$ and so  $\lbrack T,Z \rbrack = aX'$. 
        
        We also have $\lbrack T, aX' \rbrack =ac T$ and the same proof yields  that $a=0$ or $c=0$. 
      
      Up to an automorphism of $\mathcal G$, if $a \neq 0$  we can assume $a=1$,  and if  $c \neq 0$ we can  assume $c=1$.
      
      Summarizing, we have the following three posibilities concerning  the Lie algebra structure of $\mathcal G$:
      
      \begin{enumerate}
      
      \item 
      If  $a=0$ and  $c=0$, the Lie bracket relations are the following: $\lbrack Y, Z \rbrack = Z, \lbrack Y, T \rbrack = -T, \lbrack T, Z \rbrack =0$ and $\lbrack T, X' \rbrack =0$. Thus  $X'$ is  central  in $\mathcal G$. The Lie group generated by $\{ Y,Z,T  \}$  is   isomorphic to $SOL$. It then follows  that $G$ is isomorphic to the direct product 
      $\CC \times SOL$, where  $X'$ generates the center. The isotropy $I=exp(\CC Y)$ lies in $SOL$.

      \item                                                                  
      If $a=1$ and  $c=0$ the Lie bracket relations are   $\lbrack Y, Z \rbrack = Z, \lbrack Y, T \rbrack =-T, \lbrack T, Z \rbrack =X'$ and $\lbrack T, X' \rbrack =0$. The corresponding Lie
      group $G$ is isomorphic to the semi-direct product
      $\CC \ltimes Heis$, where the Lie algebra $heis$ of Heisenberg is generated by $X',T$ and $Z$.
      
      The first  factor $\CC$ is the isotropy $exp(\CC Y)$,  and its action 
      on  $heis$  is given by  $(X',Z,T) \to (X',e^t Z, e^{-t} T)$, where $X'$ is the generator of the center of $heis$. It follows that $X'$ is central  in       $\mathcal G$.
   The factor   $Heis$, intersects trivially the isotropy and hence acts freely and transitively on $G/I$.

      \item 
      For  $a=0$ and  $c=1$,  we have: $\lbrack Y, Z \rbrack = Z, \lbrack Y, T \rbrack = -T, \lbrack T, Z \rbrack =0, \lbrack T, X' \rbrack =T$      and the Lie group $G$ is a semi-direct product
      $G= \CC^2 \ltimes \CC^2$. The infinitesimal action of the first copy of $\CC^2$ (generated by $Y$ et $X'$) on the second copy of $\CC^2$ (generated by 
      $Z$ and  $T$) is given by the matrices  $\left(  \begin{array}{cc}
                                                                 1   &   0\\
                                                                 0     &  -1\\
                                                                 
                                                                 \end{array}  \right) $      and $\left(  \begin{array}{cc}
                                                                 0   &   0\\
                                                                 0     &  -1\\
                                                                 
                                                                 \end{array}  \right) $.    \end{enumerate}

\end{proof}

\subsection{The case: $H=Heis$}

Under this assumption, we will describe first the geometry of the foliation $\mathcal F$ and then we will find all algebraic  models $(G, G/I)$.

\begin{proposition}   \label{X parallel}

(i) The isotropy $I$ is unipotent.

(ii)  The $\mathcal F$-leaves are flat and $X$ is parallel along them.

\end{proposition}

\begin{proof}

The action of the isotropy $I$   on $\mathcal H / \mathcal I$  doesn't preserve any non trivial splitting. It follows that $I$ is unipotent and
$I$   is different from  the center of $H$ (which acts trivialy). This implies that any copy of $\CC^2$ transverse to the isotropy $I$ in 
$H$ acts freely and transitively on $H/I$ (they exist since $I$ is not central). This means that the $H$-leaves are flat (see the proof of Proposition~\ref{stabilizer})  and that $X$ is parallel along  them. 
\end{proof}

\begin{proposition}

$H$ is a normal subgroup of $G$.

\end{proposition}

\begin{corollaire} The $H$-foliation coincides with $\mathcal F$.
\end{corollaire}

\begin{proof}

At the Lie algebra level we show that $\mathcal H$ is an ideal in $\mathcal G$.
 Take  $A \in \mathcal G$ and let  $B$ be a local holomorphic vector field tangent to $X^{\bot}$ (recall $T\mathcal F=X^{\bot}$).
        We have to prove prove that  $\lbrack A, B \rbrack =\nabla_{A}B - \nabla_{B}A$ lies in $X^{\bot}$. 
            Note that $\nabla_{A}B \in X^{\bot}$: 
        $g(B,X)=0$ $\Longrightarrow$ $g(\nabla_{A}B,X) =-g(\nabla_{A}X,B)=0$ (because  $\nabla_{A}X = \alpha X$ by Proposition \ref{invariance}).
          On the other hand the Killing field $A$ preserves  $X$ and thus  $\nabla_{X}A = \nabla_{A}X$. As  $\nabla_{\cdot} A$ is skew-symmetric, it follows that  $g(\nabla_{B}A, X)=-
        g(B, \nabla_{X}A)=-g(B, \nabla_{A}X)=0$,  because $\nabla_{A}X = \alpha X$. The second term $\nabla_{B}A$ lies in $X^{\bot}$, and thus 
          $\lbrack A, B \rbrack \in X^{\bot}$.
          \end{proof}
          
  {\bf Algebraic  structure of $G$.}     Therefore,  $G$ is an extension of the Heisenberg group $H$. In order to describe the algebraic structure of this extension denote by $\{X', Y, Z \}$ a basis of the Lie algebra
          $\mathcal H$ of $H$, such that $Y$ is a generator of the isotropy $\mathcal I$, $X'$ is a generator of the center and $Z$ is such that: $\lbrack Y, Z \rbrack =X'$. We can assume
          that $X'$ and $Z$ generates the group of translations on the $H$-leaves.

Denote by $T$ a fourth generator of $\mathcal G$. The action of the isotropy $\CC Y$ on $\mathcal G /  \CC Y$  is such that $ad(Y)T=-Z$, which implies
$\lbrack Y, T \rbrack =-Z +  \beta Y$, for some $\beta  \in \CC$.

As the adjoint transformation of $T$ acts on $\mathcal H$ preserving the center of $\mathcal H$ it follows that: $\lbrack T, X' \rbrack =cX'$, for some constant $c \in \CC$.

We have the following

\begin{proposition}  \label{X' and X}

(i)  There exists a $H$-invariant holomorphic function on $G/I$ such that $X'=fX$ ($f$ is only locally defined on $M$ and constant on the leaves of $\mathcal F$).

(ii) $X$  is Killing (and $f$ is constant) if and only if $c=0$.

(iii) In the basis $\{ X', Z, Y \}$ of $\mathcal H$ the action of $T$ is given by $ad(T)=\left(  \begin{array}{ccc}
                                                                 c   &   m &    0\\
                                                                 0     &  c + \beta  &    1\\
                                                                 0     &  k  &  - \beta \\
                                                                 \end{array}  \right),$  with $m,k \in \CC$.
                                                                 
 (iv) If $c=0$ and $k+ \beta^2=0$, then  $g$ is flat.
\end{proposition}

\begin{proof}

(i) As $X'$ is in the center of $\mathcal H$ and $\lbrack T, X' \rbrack =cX'$,   the direction $\CC X'$ is $\mathcal G$-invariant. But in the case of unipotent isotropy
the only direction in $TM$ which is $\mathcal G$-invariant is $\CC X$. Hence $X'=f \cdot X$, for some local holomorphic function $f$ on $M$.

    Moreover, the action of $\mathcal H$ is transitive on each leaf of $\mathcal F$ and preserves $X'$ and $X$. It follows that $f$ is constant on the leaves of $\mathcal F$.
    
 (ii) As $\mathcal G$ preserves $X$, the vector field $X$ is Killing if and only if it represents   a non trivial element in  the center of $\mathcal G$. It follows that $X$ is Killing
 if and only if it 
 is a multiple of $X'$ and $X'$ is in the center of $\mathcal G$.  Equivalently, $X'$ is a central element of $\mathcal G$ if and only if $c=0$.
 
 (iii)  We apply the Jacobi relation to the vector fields  $Y,T$ and $Z$ to verify  that $ad(T)$ is a derivation if and only if $ad(T)Z$ is of the form $mX'+(c + \beta) Z+kY$, for $m,k \in \CC$.
 
 (iv) If $c=0$ and $k + \beta^2=0$, then  the vector fields  $X',Z  - \beta Y$ and $T$ generate a Lie algebra isomorphic to the Heisenberg algebra $heis$, which acts freely and transitively on $M$. The center of this algebra is generated by $X'$, which is collinear to $X$ and hence isotropic. Then,  $g$ is locally modelled on a left invariant holomorphic Riemannian metric on the Heisenberg group which gives to the center of $heis$ the norm $0$. These metrics are known to be flat~\cite{Ra}.
\end{proof}

\section{Unipotent isotropy}

In this section  we deal to  the case where the isotropy $I$ is unipotent (and $G$ is $4$-dimensional and solvable). Then,  Propositions~\ref{semi-simple classification} and~\ref{X parallel} show that $H$ is isomorphic to the Heisenberg group. The  section is devoted to the proof
of the following:

\begin{proposition}   \label{half principal}
Up to a finite unramified cover, $M$ is a quotient of $SOL$ by some lattice (and $c \neq 0$).
\end{proposition}

\subsection{Completeness.}

 Each leaf $F$ of the $H$-foliation is a surface, on which the  restriction of the vector field $X$   is an (isotropic) Killing field for the $(H,H/I)$-structure (of the leaf). The vector field $X$ generates the kernel
$\mathcal D$  of the restriction of the metric $g$ to the $F$. Furthermore,  $g$  determines  a transverse holomorphic  Riemannian structure  on the   foliation $\mathcal D$ (restricted to  $F$), i.e. 
a $(\CC, \CC)$-structure. For the basic facts concerning the study of foliations having a tranverse $(G,G/I)$-structure one can see~\cite{Mo}.
 
\begin{lemma} \label{H-completness} 
 (i)  The leaf $F$ is 
 $(H,H/I)$-complete, that  is,  the developing  map $\tilde{F} \to H/I$, on the universal
 cover,  is a diffeomorphism. 
 
 (ii) The $(G,G/I)$-structure of $M$ is complete.

\end{lemma}

\begin{corollaire} The holonomy $\Gamma$ acts properly on $G/I$.

\end{corollaire}

\begin{proof}  

(i) The $(H, H/I)$-structure on $F$ is a combination of the Killing filed $X$ and its  transverse 
$(\CC, \CC)$-structure.  One directly sees, since $X$ is complete (by compactness of $M$),
that it suffices to prove completeness of the  transverse $(\CC, \CC)$-structure, i.e. completeness of the  1-dimensional holomorphic Riemannian metric induced on  the quotient of $F$ by $X$ (or say, to prevent  any  pathology, the quotient of $\tilde{F}$ by $\tilde{X}$, where $\tilde{X}$ is the pull-back of $X$ on $\tilde{F}$).

We will show that for any complex $a$, there 
is a {\it complete} vector field $V_a$ on $F$ with (constant) $g$-norm $a$.
This would prove completeness, since such $V_a$ come from translation vector fields
on $\CC$, and hence the $V_a$'s commute, and they define a (complete) action of 
$\RR^2$, and thus the leaf is homogeneous. This action commute with the developing map, which must be diffeomorphic.

In order to check existence of the complete vector fields $V_a$, we come back to our ambient compact manifold $M$ and consider the space of vectors tangent to the $H$-foliation 
and having a norm $a$. For $a = 0$, this space is the vector bundle $\CC X$ which is known to have the global section $X$. For $a \neq 0$, this space is a fiber bundle over $M$, with fiber two copies of $\CC$ (endowed with a   structure of  an affine space).   Up to a double cover, this bundle 
is trivial and provides  a global vector field  of norm $a$ on $M$, and hence complete, by compactness of $M$.

 (ii) Since  $\mathcal H$ is an ideal of $\mathcal G$, the $H$-foliation  has a transverse $(\CC, \CC)$-structure, which  is complete by compactness of $M$~\cite{Mo}.   Combined with the completeness of the leaves, this proves completeness of the full $(G, G/I)$-structure.
 
   \end{proof}

We can  now prove:

\begin{lemma} \label{nilpotent holonomy} 

(i)  $\Gamma$ is  not abelian.

(ii) If $c=0$,  then $\Gamma$ is not nilpotent.
\end{lemma}

\begin{proof}

Consider $\overline{\Gamma}$ the  complex Zariski closure of $\Gamma$ in $G$.  As $\overline{\Gamma}$  has finitely many
connected components, up to a finite cover of $M$, we may assume that the complex abelian Lie group $\overline{\Gamma}$ is connected.

Let us notice that $\overline{\Gamma}$ can not be contained in $H$. Indeed, if not, we get a well defined surjectif projection map $M \to \overline{\Gamma} \backslash G/I \to H \backslash G/I$.  Since $I$ is
contained in $H$ and $H$ is  normal,  this last space is $\CC = H \backslash G$. This contradicts the compactness of $M$.

(i) Assume  by contradiction $\Gamma$ is abelian. Then $\overline{\Gamma}$ is an abelian complex Lie group on which the action of $\Gamma$ by adjoint representation is trivial.

Suppose first that the complex dimension of $\overline{\Gamma}$ is $1$. As above, we get a projection from $M$ to a double coset space $\overline{\Gamma} \backslash G /I$. Here $\overline{\Gamma}$ is a one-parameter complex group not included in $H$ and this double coset space
is diffeomorphic to $H/I$ which is not compact. We get a contradiction.

Assume now that the complex dimension of $\overline{\Gamma}$ is $>1$. Any element of $\overline{\Gamma}$ is invariant by the holonomy $\Gamma$ and  it gives a globally defined holomorphic Killing field on $M$. With our assumption, $M$ possesses at least two linearily independent
holomorphic (Killing) vector fields and we can use  Lemma~\ref{2 vector fields}. It follows that $M$ is a quotient of a 3-dimensional complex Lie group $C$, by a lattice $\Gamma$. As $\Gamma$ is supposed to be abelian, $C$ is also abelian and isomorphic to $\CC^3$. The holomorphic
Riemannian metric $g$ is left invariant on $\CC^3$ and hence it is flat. This is absurde, since the Killing Lie algebra   $\mathcal G$ of the flat model is of dimension $6$ (and not of dimension $4$).

(ii) Assume,  by contradiction, $\Gamma$ is nilpotent.  Since $\Gamma$ is not abelian, and supposed   to be
  nilpotent,  $\overline{\Gamma}$ is 3-dimensional (because the full group $G$ is not nilpotent)  and hence it is a complex Heisenberg  Lie group, and 
  its center is generated by $X'$.   
Take two linearily independent elements in the quotient of the Lie algebra of $\overline{\Gamma}$ outside its center. A straightforward  computation (modulo $\CC X'$) gives  $\lbrack T+aY+bZ, T+a'Y+b'Z \rbrack =(a-a')(Z - \beta Y)+ (b-b')(kY  + \beta Z)$, for all  $a,a',b,b' \in \CC$ and shows that the Lie bracket of two such elements can be  a multiple of $X'$ only if 
the determinant $k+ \beta^2$ of $\left(  \begin{array}{cc}
                                                                 1   &   -\beta\\
                                                                 \beta     &  k\\
                                                                 
                                                                 \end{array}  \right) $ vanish. Then Proposition~\ref{X' and X} implies that $g$ is flat: absurde.
\end{proof}

\subsection*{Sub-holonomy group  $\Delta = \Gamma \cap H$}

Let $\overline{\Delta}$ be the real Zariski closure of $\Delta$ in $H$. Denote by $\mathcal{\delta}$
the real Lie algebra of $\overline{\Delta}$, by $\mathcal \delta_{\CC}$ its complexified Lie algebra and by $\overline{\Delta}_{\CC}$ the associated complex Lie group.

Recall that $\Gamma$ acts on $G$ by adjoint representation and has to preserve $\Delta$ and hence also $\overline{\Delta}$ and  $\overline{\Delta}_{\CC}$.

\begin{proposition} (i) $\Delta$ is not trivial and acts properly on $H/I$.

(ii)  $\overline{\Delta}$ is of (real) dimension $\leq 4$.

(iii) $\overline{\Delta}_{\CC}$ is of (complex) dimension $\leq 2$.

(iv) $\Delta$ is abelian.
\end{proposition}

\begin{proof}

(i) Assume, by contradiction,  that  $\Delta$ is trivial. Then the projection of $\Delta$ on  $G/H \simeq \CC$ is injective  and $\Delta$
is abelian. Then Lemma~\ref{nilpotent holonomy} implies $g$ is flat: absurde.

Since  the $(H,H/I)$-structure of a leaf  $F$   is  complete, $\Delta$ is a discrete subgroup of $H$ acting properly on $H/I$ and the $\mathcal F$-leafs are diffeomorphic to 
 $\Delta \backslash H/I$.

(ii) As $H$ is  nilpotent, 
$\Delta$ is also a nilpotent group and by Malcev Theorem $\Delta$ is a (co-compact) lattice in its real Zariski closure $\overline{\Delta}$~\cite{Rag}. This means that $\overline{\Delta}$ acts properly on $H / I$
as well. Thus  $\overline{\Delta}$ has to intersect trivially the isotropy group $\CC Y \simeq \RR Y \oplus \RR iY$.

 It follows that $\overline{\Delta}$ is a real    Lie group of 
dimension $\leq 4$.

(iii) A one-parameter complex group $I'$ in $H$, not included in the subgroup of translations of $F$,  has a fix point at $x_{0}' \in F$:  it
coincides with  the isotropy at $x_{0}'$. As before, the isotropy at  $x_{0}'$ intersects trivially $\overline{\Delta}$. It follows that $\overline{\Delta}$ lies in the complex  Lie group
of translations,  whose  Lie algebra is $\CC X' \oplus \CC Z$. This implies $\mathcal{\delta}_{\CC}  \subset \CC X' \oplus \CC Z$ and  $\overline{\Delta}_{\CC}$ is of  dimension $\leq 2$.

(iv). We have $\Delta \subset \overline{\Delta}_{\CC}$, which is abelian by the previous point.
\end{proof}

\begin{proposition}   \label{dim 4}

The following facts are equivalent:

(i) The $\mathcal F$-leaves are compact;

(ii) $\overline{\Delta}$ is of (real) dimension $4$;

(iii) The projection of $\Gamma$ on $G/H$ has a discrete image.

In this case $M$ is biholomorphic to a holomorphic bundle over an elliptic curve 
with fiber type $F$ isomorphic to a $2$-dimensional complex torus.
\end{proposition}

\begin{proof}
 The $\mathcal F$-leaves  are diffeomorphic to $\Delta  \backslash H/I$. Since $\overline{\Delta}$ intersects trivially the isotropy, the  action of $\overline{\Delta}$ on $H/I$  is free
 and give a trivial foliation of $\Delta  \backslash H/I$ with compact leaves (diffeomorphic to $\Delta \backslash  \overline{\Delta}$). It follows that $\Delta  \backslash H/I$ is compact
 if and only if the action of  $\overline{\Delta}$ on $H/I$ is transitive which  means that the dimension of $\overline{\Delta}$ is $4$.
 
 The image of $\Gamma$ by the  projection $G \to G/H$  is the holonomy of the transverse $(\CC, \CC)$-structure of the $H$-foliation $\mathcal F$. The image of $\Gamma$ in $G/H \simeq \CC$ is  discrete if and only if  the leaves of $\mathcal F$ are compact~\cite{Mo}.
 
 In this case, the general study  of the developing  map of the  $(\CC, \CC)$-transverse structure of $\mathcal F$ shows that $M$ is a bundle over an elliptic curve  with fiber $F$~\cite{Mo}. 
 
Since the leaves $F \simeq \Delta \backslash  \overline{\Delta}$  are complex surfaces, $\overline{\Delta}$ is also a complex group: $\overline{\Delta}= \overline{\Delta}_{\CC}$.  It follows that $\overline{\Delta}_{\CC} \simeq \CC^2$ and  $F$ is diffeomorphic to $\Delta \backslash \CC^2$,  which is a complex torus.
\end{proof}

\begin{proposition}    \label{k}
If the complex dimension of $\overline{\Delta}_{\CC}$ is two, then $k=0$. It follows that at least one of the parameters $c$ and $\beta$ are $\neq 0$ (see Proposition~\ref{X' and X}).
\end{proposition}

\begin{proof}

Here we have  $\mathcal{\delta}_{\CC} =  \CC X' \oplus \CC Z$.

Take $\gamma \in \Gamma$ not included in $H$ and  decompose it as $\gamma=exp(\alpha T)h$, with $h \in H$ and $\alpha \in \CC^*$. 

The holonomy group  $\Gamma$ lies in the normalizer $N_{G}(\overline{\Delta}_{\CC})$ of  $\overline{\Delta}_{\CC}$
in $G$. The group $H$ normalize $\overline{\Delta}_{\CC}$ in $G$. We have then $exp(\alpha T) \in N_{G}(\overline{\Delta}_{\CC})$. It follows that the action of $ad(T)$ on $\mathcal G$ preserves $\CC X' \oplus \CC Z$. Since (by Proposition~\ref{X' and X}) we have
$\lbrack T, Z \rbrack =mX'+(c + \beta) Z +k Y$,  this implies $k=0$. Moreover, if $c= \beta=0$, then Proposition~\ref{X' and X} implies $g$ is flat: absurde.
\end{proof}

\begin{proposition}   \label{dim}
$\overline{\Delta}$  is of (real) dimension $4$.
\end{proposition}

\begin{proof}

Assume, by contradiction, $\overline{\Delta}$ is of dimension $< 4$.  Up to a finite cover, $\overline{\Delta}$ is supposed to be connected.

\subsection*{The case: $\overline{\Delta}$ is 1-dimensional} 

Then $\Delta$ is a discrete subgroup  (isomorphic to $\ZZ$) of a real one parameter subgroup $\overline{\Delta}$ of $H$. 

As $\ZZ$ does not  admit non trivial automorphisms other than $z \to -z$, up to index $2$, the action of $\Gamma$ on $\Delta$ is trivial. This implies
that the action of $\Gamma$ on $\overline{\Delta}$ is trivial as well,  and any  infinitesimal generator $Z'$ of $\overline{\Delta}$ is  an element of the real Lie algebra
$\mathcal G$ fixed by the holonomy. This element (seen as an element of the complex Lie algebra $\mathcal G$)  gives  a global holomorphic Killing field on $M$.

If the Killing field is a constant multiple of $X$,  then $c=0$ and $X$ is given by a central  element of $\mathcal G$.  It follows then that $\Delta$ lies in the center of $G$ and hence
in the center of $\Gamma$. As $\lbrack \Gamma, \Gamma \rbrack \subset \Delta$, the holonomy $\Gamma$ is a (two step) nilpotent group and Lemma~\ref{nilpotent holonomy}
gives a contradiction.

 Assume now the previous Killing field is   not  colinear
with $X$.  Note that  $\Gamma$ lies in the centralizer $C$  of $Z'$.
Since $Z'$ is not a multiple of $X'$,  the centralizer $C$ of $Z'$ is at most 3-dimensional. It follows that, up to a finite cover, $M$ admits a $(C,C)$-structure and  $M$  is  a quotient of $C$ by a lattice.

The Lie algebra of $C$ is generated by $Z'$, $X'$ and some element $T' \in \mathcal G$ not contained in $\mathcal H$. We can assume that $T'=T$(modulo $\mathcal H$). In the  Lie algebra of $C$,  the element $Z'$ is central, and $\lbrack T', X' \rbrack =c X'$. If $c \neq 0$, then $C \simeq \CC \times AG$, which is impossible since this
group is not unimodular and has no  lattices. 

It follows that $c=0$ and $C \simeq \CC^3$, which implies $g$ is flat: absurde.

\subsection*{The case: $\overline{\Delta}$ is 2-dimensional}

 The complex dimension of  $\overline{\Delta}_{\CC}$ is $1$ or $2$.

We assume first that $\overline{\Delta}_{\CC}$ is 1-dimensional. In this case $\mathcal{\delta}= \RR X'' \oplus \RR iX''$, for some  $X'' \in \mathcal G$. The adjoint action of $\Gamma$ on $\mathcal{\delta}_{\CC}= \CC X''$ is $\CC$-linear and preserves the lattice
$exp^{-1}(\Delta)$. It follows that each element of $\Gamma$ acts on $\mathcal{\delta}_{\CC}$  by homotheties given by roots of unity of order at most $6$. Up to a finite covering of $M$,
the holonomy $\Gamma$ preserves $ X''$ which gives a globally defined holomorphic Killing field on $M$. We conclude then as in  the  1-dimensional case

Assume now  $\overline{\Delta}_{\CC}$ is 2-dimensional: $\mathcal{\delta}_{\CC}=\mathcal{\delta} \otimes \CC= \CC X' \oplus \CC Z$.

          We   show that  an element $\gamma \in \Gamma$, not contained $H$,  acts trivially on $\CC X'$ and on $\mathcal{\delta_{\CC}} /  \CC X'$,
 as soon as its  projection  on $G/H$ is small enough. Such elements $\gamma$  exist, since, by  Proposition~\ref{dim 4}, the image of $\Gamma$ in $G/H$ is not discrete.

 Consider $\gamma_{n}=r_{n}h_{n}$ a sequence of elements of $\Gamma$, with $h_{n} \in H$ and $r_{n} \notin H$ going to $0$ in $G/H \simeq \CC$, when $n$ goes to infinity.
We can assume $r_{n}=exp(\alpha_{n}T)$, with $\alpha_{n} \in \CC^*$ going to $0$ when $n$ goes to infinity.

If $h_{n}=exp(a_{n}X')exp(b_{n}Y)exp(c_{n}Z),$ with $a_{n},b_{n},c_{n } \in \CC$ then the adjoint action of $h_{n}$ on $\overline{\Delta}_{\CC}$  is exactly the action of $Ad(exp(b_{n}Y))$.

The action of $Ad(exp(b_{n}Y))$ on $\mathcal{\delta}_{\CC}= \CC X' \oplus \CC Z$  is given by the matrix  $\left(  \begin{array}{cc}
                                                                 1   &   b_{n} \\
                                                                 0     &  1\\
                                                                 \end{array} \right)$.

      By Proposition~\ref{X' and X}, $Ad(r_{n})=Ad(exp(\alpha_{n}     T))$ has the following matrix when acting on $\mathcal{\delta}_{\CC}= \CC X' \oplus \CC Z$:
      $\left(  \begin{array}{cc}
                                                                e^{\alpha_{n} c} &  * \\
                                                                 0     &  e^{\alpha_{n} (c + \beta)}\\
                                                                 \end{array} \right)$. The matrix of $Ad(\gamma_{n})=Ad(r_{n})Ad(h_{n})$ has    the same form.
           
           Recall now that this   action of $Ad(\gamma_{n})$ preserves $\mathcal{\delta}$ and the lattice $exp^{-1}(\Delta)$: it is conjugated to an element of $SL(2, \ZZ)$.
           It follows that,  for all $n \in \NN$, the
            previous matrix  of $Ad(\gamma_{n})$ has a determinant which equals  $1$ and a trace which is an integer. This implies that,  for $n$ large enough,
           the trace equals $2$ and  $e^{\alpha_{n}c}= e^{\alpha_{n} (c  + \beta)}=1$. It follows $c=0$ and $\beta =0$,  which   contradicts   Proposition~\ref{k}.

 \subsection*{The case:  $\overline{\Delta}$ is 3-dimensional}
        
    As in the previous case, we have  $\mathcal{\delta}_{\CC}= \CC X' \oplus \CC Z$. We can change the infinitesimal generator $X'$ of the center of $H$ and also $Z$ into $Z+ a X'$, with $a  \in \CC$, such that  either $\mathcal{\delta}= \CC X' \oplus \RR Z$,         or $\mathcal{\delta}= \RR X' \oplus \CC Z$. The previous transformation  keeps unchange the Lie bracket relations.

     Take as before a sequence  $\gamma_{n} =exp(\alpha_{n}T)h_{n}$ of elements of $\Gamma$, such that $h_{n} \in H$ and $\alpha_{n}  \in \CC^*$ converges to $0$. As before,
     the matrix of the $Ad(\gamma_{n})$-action on $\mathcal \delta_{\CC}= \CC X' \oplus \CC Z$  is of the form $\left(  \begin{array}{cc}
                                                                e^{\alpha_{n} c} &  * \\
                                                                 0     &  e^{\alpha_{n} (c  + \beta)}\\
                                                                 \end{array} \right)$.
                                                                 
       Consider  the restriction of $Ad(\gamma_{n})$  to  $\mathcal \delta$. For each $n \in \NN$, the $Ad(\gamma_{n})$-action on $\mathcal \delta$ preserves some lattice, so it is conjugated to some element in $SL(3, \ZZ)$. When $n$ goes to infinity, the three eigenvalues
       of $Ad(\gamma_{n})$ go to $1$. By discretness of $SL(3, \ZZ)$, it follows that, for $n$ large enough,  all  eigenvalues of $Ad(\gamma_{n})$ equal $1$. So, for
       $n$ large enough, $e^{\alpha_{n}c}= e^{\alpha_{n} (c  + \beta)}=1$. It follows $c=0$ and $\beta =0$,  which  contradicts 
       Proposition~\ref{k}.
        \end{proof}

 We are now able to prove Proposition~\ref{half principal}.
 
 \begin{proof}
By  Proposition~\ref{dim 4}, $M$ is a fiber bundle over an elliptic curve with fiber $F$ biholomorphic to a $2$-dimensional complex torus. We have seen that $\Delta$ is an abelian group isomorphic to $\ZZ^4$,
  $\overline{\Delta} \simeq \RR^4$ and $\overline{\Delta}_{\CC} = \CC^2$. As before, we have $\mathcal{\delta}_{\CC} =  \CC X' \oplus \CC Z$.

    By Proposition~\ref{dim 4},  the projection of $\Gamma$ on $G/H$ is a discrete subgroup. This subgroup is isomorphic to the fundamental group of the basis of our fibration,  so it is $\simeq \ZZ^2$. Take $\gamma_1$ and $\gamma_{2}$ two elements in $\Gamma$ such that their
  projections in $G/H$ span the previous $\ZZ^2$. Then any element of $\Gamma$ decomposes  as $\gamma_1^p \gamma_2^q d$, with $p,q \in \ZZ$ and $d \in \Delta$. Moreover,
  we can decompose $\gamma_{i}$ as $exp(\alpha_{i}T)h_{i}$, where $i \in \{1, 2 \}$, $h_{i} \in H$ and $\alpha_{i} \in \CC$.

  Assume by contradiction that $c=0$. Then Proposition~\ref{X' and X} implies that the action of $Ad(T)$ on the quotient $\mathcal H /  \mathcal I$ is of (complex) determinant $1$.
  Hence  the determinant of the action of $Ad(\gamma_{i})$ on $\mathcal{\delta}_{\CC}$ equals $1$.
  
  On the other hand the eigenvalues of $Ad(\gamma_i)$ are $1$ and $e^{\alpha_{i} \beta}$ (see the proof of the case $2$ in Proposition~\ref{dim}). It follows that $e^{\alpha_{i} \beta}=1$, for $i \in \{1,2 \}$. This implies 
  $\alpha_{i} \beta =2i \pi k_{i}$, where $k_{i} \in \ZZ$. Since $\alpha_{i}$ are $\ZZ$-independent, we have $\beta =0$. As before, this implies $n=0$ and $g$ is flat: absurde.

  It follows that  $c \neq 0$.

  We prove  that there exists a  basis of $\mathcal{\delta}_{\CC}$ in respect of which  the actions of $Ad(\gamma_{1})$ and $Ad(\gamma_2)$ are  (both) diagonal. Recall that $\CC X'$ is stable by the adjoint representation of $G$ and, in particular, by $Ad(\gamma_{1})$ and by
 $Ad(\gamma_{2})$. Denote $\lambda_i$ the corresponding eigenvalue 
  of the restriction of $Ad(\gamma_i)$ to $\mathcal{\delta}_{\CC}$, $i \in \{1,2 \}$.  We prove by contradiction that either the modulus of $\lambda_1$ or the modulus of $\lambda_2$  is $\neq 1$. Indeed, if not the  modulus of the "quotient" $f$ of $X'$ over $X$ (see Proposition~\ref{X' and X}) is preserved
  by the projection of $\Gamma$ on $G/H$ (which coincides with the holonomy of the transversal structure of the $H$-foliation). This means $|f|$ is globally defined on $M$. As $M$ is compact and $f$ is holomorphic,  the maximum principle implies $f$ is constant and,  by Proposition~\ref{X' and X},  we have $c=0$, which contradicts our assumption.
  
  Assume now that the modulus of $\lambda_1$ is $\neq 1$. As $Ad(\gamma_1)$ acts on $\mathcal{\delta_{\CC}}$ preserving a lattice,  this action is unimodular. It follows that the action of $Ad(\gamma_1)$ on $\mathcal{\delta}_{\CC}$ has distinct eigenvalues, and  so it is diagonalizable over $\CC$. Since $\gamma_1$ and
  $\gamma_2$ commutes (modulo $\Delta$) and the action of $\Delta$  on $\mathcal{\delta}_{\CC}$ is trivial, then $Ad(\gamma_1)$ and $Ad(\gamma_2)$ commutes in restriction to $\mathcal{\delta}_{\CC}$. It follows that the two eigenvectors of $Ad(\gamma_1)$ are invariant  by $Ad(\gamma_2)$ as well.
  Consequentely the two eigenvectors of $Ad(\gamma_1)$ are $\Gamma$-invariant.
   The holonomy group $\Gamma$ lies in a  subgroup of $G$ for which the adjoint action on $\mathcal{\delta}_{\CC}$ preserves a  non trivial splitting.
  
  Take $T' \in \mathcal G$ such that $\gamma_1=exp(T')$. We have  proved that $\Gamma$ lies in the $3$-dimensional (solvable) complex Lie group $C$ generated by $\CC T'$ and $\mathcal{\delta}_{\CC}$.  Thus,  the manifold $M$ possesses a $(C,C)$-structure
 and $M$ is a quotient of $C$ by a lattice (so $C$ is unimodular). Since $c \neq 0$,  the only compatible Lie group structure is $SOL$ and so, up to a finite cover,  $M$ is a quotient of $SOL$ by some lattice.
  \end{proof}

       \section{Semi-simple isotropy}   \label{semi-simple isotropy}
       
       \subsection{Solvable Killing algebra} We study separately the  two possible models we get in Proposition~\ref{semi-simple classification}. We prove the following:
       
       \begin{proposition}
      Up to a finite unramified cover, $M$ is a quotient of the Heisenberg group by a lattice ($G$ is isomorphic to $\CC \ltimes Heis$).
   \end{proposition}
      
      Together with Proposition~\ref{half principal} this will prove part (ii) ($\mathcal G$ solvable) of the main  Theorem \ref{result}.

{\bf The Case  $G= \CC \times SOL.$}
    
Recall  the Lie algebra of  $SOL$ is generated by  $\{ Z,T,Y \}$,  with the Lie bracket relations $\lbrack Y, Z \rbrack = Z, \lbrack Y, T \rbrack = -T$ and $\lbrack T, Z \rbrack =0$. The center of $\mathcal G$  is generated by $X'$ and the 3-dimensional abelian  Lie algebra generated by
  $\{ X',Z,T  \}$ acts freely and transitively on $G/I$. The holomorphic Riemannian metric $g$ is locally identified with a translation-invariant holomorphic Riemannian metric on $\CC^3$.
  Consequently $g$ is flat, which is impossible.

{\bf The case $G= \CC \ltimes Heis.$}                                                          
      
      Recall that the Lie algebra of $Heis$ is generated by the central element $X'$ and by $Z,T$ such that $\lbrack Z, T \rbrack =X'$. We have seen   that $X'$ is    fixed by the isotropy $I$
      and $Z$ and $T$ are the two isotropic directions expanded  and contracted     by $I$.
                                                       
  Here  $X'$ generates the  global Killing   field $X$ of constant norm  equal to $1$ fixed by the isotropy. Denote $\phi^t$, where $t \in \CC$,  the holomorphic flow of $X$.
  The flow $\phi^t$ preserves the  orthogonal  distribution $X^{\bot}$. This distribution has dimension 2 and it is non-degenerate in respect to $g$. Thus 
  $X^{\bot}$  has exactly  two isotropic line fields which are locally generated by $Z$ and $T$. They are naturally preserved 
  by $\phi^t$. Since $\lbrack Z, T \rbrack \neq 0$, the distribution $X^{\bot}$ is not integrable.
  
 We will say that $X$ is {\it equicontinuous}  if  $\phi^t$ is. This means by definition that the closure $K$ of $\phi^t$ in the group of homeomorphisms of $M$
   is a compact group. In this case $K$ will be an abelian  compact complex Lie group (a complex torus) acting on $M$ and preserving $g$.
   
   Assume first that $X$ is equicontinuous. If $K$ has complex dimension $>1$, the fundamental fields of the action of $K$ on $M$ give    
   at least two linearily independent global holomorphic vector fields on $M$ and Lemma~\ref{2 vector fields} applies.   So the centralizer $C$ of $K$ in $G$ 
    acts transitively on $M$,  such that $M$ is quotient of $C$ by a lattice. The subgroup $C$ of $G$ is unimodular and  has a center which is at least 1-dimensional. It follows that
    $C$ is isomorphic to $Heis$.

    Now consider the case where $K$ is a 1-dimensional complex torus.
     The quotient of $M$ by the action of $K$ is a compact complex surface $S$ which inherits  a flat holomorphic Riemannian metric. Indeed, $G/exp(\CC X) \simeq SOL$ and $S$ is easily seen to be locally modelled on $(SOL, SOL/I')$, where $SOL \simeq \CC \ltimes \CC^2$
    with the action of $\CC$ on $\CC^2$ given by the complex one-parameter group $I'=\left(  \begin{array}{cc}
                                                                 e^t   &   0\\
                                                                 0     &  e^{-t}\\  \end{array}   \right).$

   Up to a finite unramified  cover,  this surface   is  a $2$-dimensional complex torus $T^2$ with a flat holomorphic Riemannian metric (see  Theorem 4.3 in \cite{Dum}).
   Consequently, up to a finite unramified  cover,  $M$ is a principal bundle of elliptic curves over  a complex torus and the projection of the holonomy $\Gamma$ on $G/exp(\CC X) \simeq SOL$ lies in the subgroup of translations $\CC^2$. It follows that the holonomy $\Gamma$ lies in a complex Lie group $C$ 
   of dimension $3$ which is a central extension of $\CC^2$ by $\CC$
  (isomorphic to $Heis$)  and which acts freely and transitively on $G/I$. Up to a finite unramified cover,  $M$ is biholomorphic to a quotient of $Heis$ by a lattice. 
   
   It remains to settle the case where $X$ is non-equicontinuous, for which we prove:
   
   \begin{proposition}     \label{Anosov}
   If the flow $\phi^t$ is non-equicontinuous,  then it  is Anosov.
   \end{proposition}
   
   \begin{proof} 
   
  By passing,  if necessary,  to a finite cover,  we may assume that the two isotropic directions of $X^{\bot}$ are directed by two smooth vector fields $T_{1}$ and $T_{2}$. The $\phi^t$-invariance of these isotropic directions shows  that $D_{x} \phi^t(T_{1}(x))=a(x,t)T_{1}(\phi^t(x))$ and $D_{x} \phi^t(T_{2}(x))=b(x,t)T_{2}(\phi^t(x))$, for any $x \in M$ and $t \in \CC$;
  $a$ and $b$ being some smooth complex valued functions on $M \times \CC$. By the volume preserving property $a(x,t)b(x,t)=1$.

  We now prove that for any $x$, the orbit  $\{D_{x} \phi^t(T_{1}(x)), t \in \CC \}$ is not bounded in $TM$. Assume, by contradiction,  that
  the modulus of the function $a$ is upper bounded.  If the modulus of $a(x,t)$ stays  $\geq a' >0$ for a sequence $t_{n}$ tending to $+ \infty$
  or $- \infty$,  then $D_{x} \phi^{t_{n}}$ is equicontinuous and so by a simple result of~\cite{Zeg} the flow itself is equicontinuous, which contradicts our hypothesis. It then follows
  that $a(x,t) \to 0$, when $t \to + \infty$ or $t \to - \infty$. Thus (by continuity of $a$) there are two sequences $t_{n}$ and $t_{n'}$ tending to $+ \infty$, such that $a(x,-t_{n})=a(x,t'_{n})$.
  By the cocycle property of $a$, applied to $x_{n}= \phi^{-t_{n}}(x)$, we get: $a(x_{n}, t'_{n} +t_{n})=a(x,t'_{n})a(x_{n},t_{n})$. But $a(x_{n},t_{n})a(x,-t_{n})=1$, and hence $a(x_{n},t_{n}+t'_{n})=1$. Hence $b(x_{n},t_{n}+t'_{n})=1$, and consequently $D_{x_{n}} \phi^{t_{n} +t'_{n}}$ is equicontinuous. Since $t_{n} +t'_{n}$ tends to $+ \infty $,  Proposition
  3.2 of \cite{Zeg} implies  then that $\phi^t$ is equicontinuous which contradicts our assumption. 
  
  In the same way,  the modulus of $b$ is unbounded and hence the orbit
  of any non zero vector in $X^{\bot}$  under  the action of $D\phi^t$ is not bounded. This means, by definition,  that $\phi^t$ is quasi-Anosov and by an easy  case of the main Theorem in \cite{Man}  this implies that $\phi^t$ is a holomorphic Anosov flow in  Ghys's  sense~\cite{Ghys}.
   \end{proof}
   
   A simple  case of the classification of  holomorphic Anosov flows on compact complex $3$-manifolds~\cite{Ghys} shows that $\phi^t$ preserves some holomorphic Riemannian metric $q$ of constant sectionnal curvature. As $X^{\bot}$ is not integrable, $q$ is necessarily   of non-zero constant sectionnal curvature~\cite{Ghys}.  By Theorem~\ref{Dumitrescu}, the intersection $\mathcal G'$ of the Killing Lie algebra of $g$ and the Killing Lie algebra of $q$ acts transitively on $M$. This implies that the Heisenberg
   algebra is contained in the Killing Lie algebra $sl(2, \CC) \oplus sl(2, \CC)$ of $q$. This is absurde, and therefore,  $X$ is equicontinuous.
   
   \subsection{Semi-simple Killing algebra}
   
   Here $G=\CC \times SL(2, \CC)$ and $I=\left(  \begin{array}{cc}
                                                                 e^t  &   0\\
                                                                 0     &  e^{-t}\\
                                                                 
                                                                 \end{array}  \right) \subset SL(2,\CC).$

         We show the following 
         
         \begin{proposition}  
         There are no compact manifolds locally modelled on $(G,G/I)$.
         \end{proposition}
         
         This will complete  the proof of the main Theorem~\ref{result}.
         
         \begin{proof}
         The factor $\CC$ of $G$ is generated by the flow of the Killing vector field $X$. 
         
         Assume first that $X$ is equicontinuous and consider the complex Lie group $K$ which is the closure
         of the flow of $X$ in the group of homeomorphism of $M$. We have seen that if the complex dimension       of $K$ is $>1$ then,  Lemma~\ref{2 vector fields} implies that there exists a 3-dimensional complex subgroup  $C$ in $G$ which acts freely and transitively on $M$ and $M$ identifies with a quotient of $C$ by some lattice. This is impossible because the only 3-dimensional
         subgroups of $G$ which act freely on $M$ are isomorphic to $\CC \times AG$ and they do not  have lattices (they are not unimodular).
         
         If $K$ has dimension 1  the quotient of $M$ by $K$ is a complex compact surface locally modeled on $(SL(2,\CC), SL(2, \CC) /I)$. This compact surface possesses a holomorphic Riemannian metric of non-zero constant sectionnal curvature. But, by Theorem 4.3 in \cite{Dum}, all holomorphic Riemannian metrics on compact complex surfaces are  flat, which leads to   a contradiction.
         
         Consider now the case where $X$ is non-equicontinuous. The proof of Proposition~\ref{Anosov} implies that $X$ is an Anosov flow with stable and      instable directions
         given by the isotropic directions of $X^{\bot}$. Here the holomorphic plane field $X^{\bot}$         is integrable because it is tangent to the orbits of $sl(2, \CC)$-action. In this situation Ghys's  classification~\cite{Ghys} shows that, up to a finite cover,  $M$ is biholomorphic to a holomorphic suspension   (given by the flow of $X$)  of a complex hyperbolic
         linear automorphism of a complex torus $T^2$. In particular, the orbits of $sl(2,\CC)$ are 2-dimensional complex torii locally modelled on $(SL(2,\CC), SL(2, \CC) /I)$. We get the same contradiction as before.
        \end{proof}

${}$ 

\end{document}